\documentclass[10pt]{amsart}
\usepackage{amsmath}
\usepackage{amssymb}
\usepackage{enumerate}
\usepackage{amsbsy}
\usepackage{amsfonts}
\usepackage{color}
\usepackage{upgreek}
\usepackage{parskip}

\newtheorem{thm}{Theorem}[section]
\newtheorem{lem}[thm]{Lemma}
\newtheorem{prop}[thm]{Propsition}

\numberwithin{equation}{section}

\begin{document}

	\title[Hartree-type nonlinear Dirac equation]{Scattering and non-scattering of the Hartree-type nonlinear Dirac system at critical regularity}
	
	\author[Y. Cho]{Yonggeun Cho}
	\address{ Department of Mathematics, and Institute of Pure and Applied Mathematics, Jeonbuk National University, Jeonju 561-756, Republic of Korea}
	\email{changocho@jbnu.ac.kr}
	
    \author[S. Hong]{Seokchang Hong}
    \address{Department of Mathematical Sciences, Seoul National University, Seoul 08826, Republic of Korea}
    \email{seokchangh11@snu.ac.kr}

 \author[K. Lee]{Kiyeon Lee}
 \address{ Department of Mathematics, Jeonbuk National University, Jeonju 561-756, Republic of Korea}
 \email{leeky@jbnu.ac.kr}

	\thanks{2020 {\it Mathematics Subject Classification.} 35Q55, 35Q40.}
	\thanks{{\it Key words and phrases.} Dirac equation, global well-posedness, scattering, Yukawa potential, Coulomb potential, null structure, angular regularity, $U^p-V^p$ space.}
	
	\begin{abstract}
	We consider Cauchy problem of the Hartree-type nonlinear Dirac equation with potentials given by $V_b(x) = \frac1{4\pi}\frac{e^{-b|x|}}{|x|}\, (b \ge 0)$. In previous works, a standard argument is to utilise null form estimates in order to prove global well-posedness for $H^s$-data, $s>0$. However, the null structure inside the equations is not enough to attain the critical regularity. We impose an extra regularity assumption with respect to the angular variable. Firstly, we prove global well-posedness and scattering of Dirac equations with Hartree-type nonlinearity for $b>0$ for small $L^2_x$-data with additional angular regularity. We also show that only small amount of angular regularity is required to obtain global existence of solutions. Secondly, we obtain non-scattering result for a certain class of solutions with the Coulomb potential $b=0$.
	\end{abstract}

		\maketitle

\section{Introduction}
We are concerned with global well-posedness and scattering of the $(1+3)$-dimensional Dirac equation with Hartree-type nonlinearity for Yukawa and Coulomb potentials.
The main equation is given by
\begin{align}\label{main-eq}
\left\{\begin{aligned}
-i\gamma^\mu\partial_\mu\psi + m\psi &= (V_b*(\psi^\dagger\gamma^0\psi))\psi, \\
 \psi(0,\cdot)&:=\psi_0\in L^2_x(\mathbb R^3).
\end{aligned}
\right.
\end{align}
The potential $V_b $ is the spatial function
\begin{align}\label{yukawa}
V_b(x) = \frac1{4\pi}\frac{e^{-b|x|}}{|x|} \;\;(b \ge 0),
\end{align}
which is called Yukawa (Coulomb) potential if $b > 0\; ( b = 0)$.  

We start with the basic notation. Throughout this paper, we denote points by $(x^\mu)$, $\mu=0,1,2,3$ in the Minkowski space $(\mathbb R^{1+3},\mathbf m)$, where $\mathbf m$ is the metric given by $\mathbf m=\textrm{diag}(-1,+1,+1,+1)$. The partial derivatives with respect to $x^\mu$ is written by $\partial_\mu$. We shall use the notation $t=x^0$ for time variable, and $x=(x^1,x^2,x^3)$ for spatial variable. Then we write $\partial_0=\partial_t$ and $\nabla=(\partial_1,\partial_2,\partial_3)$. The unknown spinor field $\psi$ is written as a column vector in $\mathbb C^4$ and $m>0$ is a mass constant. We also denote the complex conjugate of the  transpose $\psi^t$ by $\psi^\dagger$. We use Roman indices $j,k=1,2,3$ and Greek indices $\mu,\nu=0,1,2,3$ and adopt the Einstein summation notation, i.e., any repeated indices mean the summation over described range. Thus we write $\gamma^\mu\partial_\mu = \gamma^0\partial_t+\sum_{j=1}^3\gamma^j\partial_j$, where $\gamma^\mu$, $\mu=0,1,2,3$ are the Dirac gamma matrices given by
\begin{align}\label{gm-matrix}
\gamma^0 = \begin{bmatrix} I_{2\times2} & 0 \\ 0 & -I_{2\times2} \end{bmatrix},\quad \gamma^j = \begin{bmatrix} 0 & \sigma^j \\ -\sigma^j & 0 \end{bmatrix}	
\end{align}
 with Pauli matrices:
 \begin{align}\label{pauli}
 \sigma^1 = \begin{bmatrix} 0 & 1 \\ 1 & 0 \end{bmatrix}, \quad \sigma^2 = \begin{bmatrix} 0 & -i \\ i & 0 \end{bmatrix}, \quad \sigma^3 = \begin{bmatrix} 1 & 0 \\ 0 & -1 \end{bmatrix}.
 \end{align}

One may observe that the equation \eqref{main-eq} can be derived by uncoupling the Dirac-Klein-Gordon system
\begin{align}\label{dkg}
\begin{aligned}
(-i\alpha^\mu\partial_\mu+m\beta)\psi = \phi\beta\psi, \\
(\partial_t^2-\Delta+M^2)\phi = \psi^\dagger\beta\psi.
\end{aligned}
\end{align}
In fact, we let a scalar field $\phi$ be a standing wave, i.e., $\phi(t,x)=e^{i\lambda t}f(x)$ with $M \ge |\lambda|$. Then the Klein-Gordon part of \eqref{dkg} becomes
\begin{align}\label{kg}
(- \Delta  + M^2 -\lambda^2)\phi = \psi^\dagger\beta\psi.
\end{align}
Then one easily shows that the solutions of \eqref{kg} are given by
\begin{align}\label{kg-sol}
\phi =  V_b *(\psi^\dagger\beta\psi)
\end{align}
with $b = \sqrt{M^2-\lambda^2}$. We put \eqref{kg-sol} into the Dirac part of \eqref{dkg} and then a spinor field $\psi$ gives the desired equation. (See also \cite{tes1,cyang}.) One may also replace a quadratic term $\psi^\dagger\beta\psi$ by $|\psi|^2$. Indeed, the equation \eqref{dirac-eq} with $V_0(x) = |x|^{-1}$ was derived by the authors of \cite{changla}, by uncoupling the Maxwell-Dirac system under the assumption of vanishing magnetic field with the quadratic term $|\psi|^2$.

The $L^2_x$-norm of the solutions to the system \eqref{main-eq} is conserved:
\begin{align}\label{m-conserv}
\int_{\mathbb R^3}|\psi(t,x)|^2\,dx = \int_{\mathbb R^3}|\psi_0(x)|^2\,dx.
\end{align}

When $m = b = 0$, the equation \eqref{dirac-eq} is invariant under the scaling:
$$
\psi(t,x) \mapsto \psi_\lambda(t,x) = \lambda^\frac32 \psi(\lambda t,\lambda x),
$$
for fixed $\lambda>0$. Thus the system \eqref{dirac-eq} is essentially $L^2$-critical.
In this paper we exclusively consider the massive case $(m > 0)$ and by scaling, we set $m=1$ hereafter.

We introduce the notation $\alpha^\mu$ and $\beta$ as follows:
\begin{align}\label{albe}
\alpha^j = \gamma^0\gamma^j, \quad \beta = \gamma^0.
\end{align}
To study the initial value problem of Dirac equations, we shall follow the standard approach as in \cite{candyherr,tes,cyang}. We define the Dirac projection operators $\Pi_\pm$ as a Fourier multiplier with the symbol
\begin{align}\label{proj}
\Pi_\pm(\xi) = \frac12\left(I_{4\times4}\pm\frac{\alpha^j\xi_j+\beta}{\langle \xi\rangle}\right),	
\end{align}
where $\langle\xi\rangle=(1+|\xi|^2)^\frac12$. A simple computation gives the following properties:
$
\Pi_\pm\Pi_\pm = \Pi_\pm,\ \Pi_\pm\Pi_\mp = 0.
$ Throughout this paper, we will use the notation $\psi_\pm := \Pi_\pm\psi$. Then we can decompose $\psi = \psi_++\psi_-$. Finally, using notation \eqref{albe} and projection \eqref{proj}, our system \eqref{main-eq} is rewritten as
\begin{align}\label{dirac-eq}
(-i\partial_t\pm \langle \nabla\rangle)\psi_\pm = \Pi_\pm[(V_b*(\psi^\dagger\beta\psi))\beta\psi],\;\;\psi_\pm(0) := \psi_{0, \pm}.	
\end{align}
We say that the solution $\psi$ scatters to a free solution in a Hilbert space $\mathcal H$ if there exist $\psi_{\pm }^\ell := e^{\mp it\langle \nabla\rangle}\varphi_{\pm } \;(\varphi_{\pm} \in \mathcal H)$ such that
\begin{align*}
\|\psi_\pm(t) - \psi_{\pm}^\ell(t)\|_{\mathcal H} \to 0\;\;\mbox{as}\;\;t \to \pm \infty,
\end{align*}
or equivalently,
\begin{align*}
\|\psi(t) - \psi^\ell(t)\|_{\mathcal H} \to 0\;\;\mbox{as}\;\;t \to \pm \infty,
\end{align*}
where $\psi^\ell = \psi_+^\ell + \psi_-^\ell$.

Recently, the Cauchy problem of Dirac equations has been extensively studied.
For instance, see \cite{bejeherr1} for result on the cubic-nonlinear Dirac equations and reference therein.
%
Dirac equations coupled with several fields have been also well-studied. For example, we refer the readers to \cite{danfoselb1,wang,candyherr,candyherr1} for the study on the Maxwell-Dirac and Dirac-Klein-Gordon systems.

We shall mention a few selected results on some related equations. Firstly, we present the boson star equation (or semi-relativistic equation) with Hartree-type nonlinearity:
\begin{align}\label{bo-star}
	(-i\partial_t+\langle \nabla \rangle)u = \left(V_b * |u|^2\right)u.
\end{align}
In the works of Lenzmann and Cho-Ozawa \cite{lenzmann, chooz} the well-posedness for $b \ge 0$ was proved in $H^s$-data, $s > \frac12 - \varepsilon$, and it was improved later to $s > \frac14$ by Herr-Lenzmann \cite{herrlenz} when $b = 0$. The linear scattering does not occur when $b = 0$ \cite{chooz}. Instead, a modified scattering of \eqref{bo-star} can occur in case when $b = 0$. For this see \cite{pusateri}. On the other hand, the linear scattering problem was handled in \cite{herrtesf} for $b > 0$ and $s > 0$.

\subsection*{Scattering for Yukawa potential}
  Now we pay attention to the equation \eqref{dirac-eq} with the Yukawa potential. The Cauchy problem of the system \eqref{dirac-eq} was studied by A. Tesfahun \cite{tes,tes1} and C. Yang \cite{cyang} independently. The authors of \cite{tes1,cyang} utilise the null structure and bilinear estimates to prove global well-posedness and scattering for $H^s$-data, $s>0$. However, the global well-posedness is still open at the critical regularity.

In this paper we establish global well-posedness and scattering of solutions to the system \eqref{dirac-eq} for small data in the scaling critical Sobolev space which has extra weighted regularity in the angular variables. To be more precise, we let $\Omega_{ij}=x_i\partial_j-x_j\partial_i$ be the infinitesimal generators of the rotations on $\mathbb R^3$ and let $\Delta_{\mathbb S^2} = \sum_{1 \le i < j \le 3}\Omega_{ij}^2$ be the Laplace-Beltrami operator on the unit sphere $\mathbb S^2\subset\mathbb R^3$. Then we can define the fractional power of angular derivative by $\langle\Omega\rangle^\sigma=(1-\Delta_{\mathbb S^2})^\frac\sigma2$, which will be treated concretely below, and define angularly regular space $L_x^{2, \sigma}$ space by $\langle \Omega \rangle^{-\sigma}L_x^2$ and its norm by $\|f\|_{L_x^{2, \sigma}} := \|\langle \Omega \rangle^\sigma f\|_{X_x^2}$.
Now we state the main theorem:
\begin{thm}\label{gwp}
Let $\sigma > 0$. Then there exists $\delta>0$ such that for initial data  $\|\psi_0\|_{L^{2, \sigma}_x(\mathbb R^3)}\le \delta$, the Cauchy problem \eqref{dirac-eq} is globally well-posed and solutions $\psi$ scatter in $L_x^{2, \sigma}$ to free solutions as $t\rightarrow\pm\infty$.
\end{thm}

The main improvement of Theorem \ref{gwp} is to attain the critical regularity. Motivated by the work of \cite{wang,candyherr}, we exploit an additional angular regularity. Furthermore, we observe that only a small amount of regularity in the angular variables is required to prove global well-posedness.

\subsection*{Strategy of proof of Theorem \ref{gwp}}
We discuss the key ideas of the proof of Theorem \ref{gwp}. 
Because of huge amount of notations to be used in the rest of this paper, we would like to elucidate the main scheme and motivation here for convenience to the readers.
The main approach is to construct the Picard's iterate, which is convergent in the adapted function spaces. Thus the crucial part of the proof is the following multilinear estimates: (See also Proposition \ref{main-est}.)
$$
\left\| \int_0^t e^{\mp i(t-t')\langle\nabla\rangle} [V_b*(\varphi^\dagger\beta\phi)\beta\psi](t')\,dt'\right\|_{F^\sigma_\pm}  \lesssim \|\varphi\|_{F^\sigma_{\pm_1}} \|\phi\|_{F^\sigma_{\pm_2}}	\|\psi\|_{F^\sigma_{\pm_3}}.
$$
Here $F^\sigma_\pm$ is the adapted function space which will be defined in Section 3. Roughly speaking, the space $F^\sigma_\pm$ consists of $V_\pm^2$ space equipped with angular regularity and hence by making use of duality (Lemma \ref{lem-v-dual}) we will study quartilinear estimates.

We deal with all possible frequency interactions such as High$\times$High and Low$\times$High interactions with low-modulation and high-modulation regimes. Since we have four-input frequencies, it seems to require repetitive work. Fortunately, by H\"{o}lder's inequality and symmetry between two spinor fields, the problem of the quartilinear estimates can be reduced to frequency-localised bilinear expressions such as
$$
\|P_{\lambda_0}[(P_{\lambda_1}\varphi)^\dagger\beta(P_{\lambda_2}\phi)]\|_{L^2_tL^2_x},
$$
 except for the case that high-modulation is bigger than the highest-input-frequency. See also Proposition \ref{bi-est}. When the modulation is bigger than the highest frequency, the situation is rather easier than other cases. Indeed, we simply use boundedness in high modulation regime \eqref{bdd-high-mod} and  $L^2$-bilinear estimates shown in \cite{cyang}.
\subsubsection*{High frequency - Low modulation}
We consider the case that the modulation $d$ is less than the lowest-input-frequency. There is nothing new ingredients to obtain the required bilinear estimates. Indeed, the space-time Strichartz estimates and the null structure between two input-spinor fields will play a crucial role in the low-modulation regime. However, the localised $L^4$-Strichartz estimate (see \cite{choozxia} for instance) gives
$$
\|e^{\mp it\langle\nabla\rangle}P_\lambda f\|_{L^4_tL^4_x} \lesssim \lambda^\frac12 \|P_\lambda f\|_{L^2_x},
$$
which would be too big and troublesome in the summation. To avoid this problem, we apply the almost orthogonal decomposition by cubes with smaller size $\mu\le\lambda$ . Then we have the improved estimates such as
$$
\|e^{\mp it\langle\nabla\rangle}P_qP_\lambda f\|_{L^4_tL^4_x} \lesssim (\mu\lambda)^\frac14 \|P_qP_\lambda f\|_{L^2_x},
$$
where $P_q$ is the cube localisation operator. Even though we gain factor $\left(\frac\mu\lambda\right)^\frac14$, it pays for more work, i.e., we need to take square-summation by cubes to recover the $\|P_\lambda f\|_{L^2_x}$ term. This step would cause some loss in a certain estimate. Here Lemma \ref{square-sum} assures that such loss can be absorbed elsewhere.

We would like to mention that the most delicate interaction is the High$\times$High frequency interaction. This is why the low regularity problem becomes more difficult as the spatial dimension decreases, i.e., the High$\times$High interactions in the nonlinearity grows seriously. Thus such interactions would be the main obstacle in the improvement of the previous results \cite{tes1,cyang}. At this point, we remark that the (infinitesimal) rotation generators $\Omega_{ij}=x_i\partial_j-x_j\partial_i$ relax such delicate interactions. Furthermore, since the optimality of available range of the Strichartz estimates is given by the Knapp-type counterexample, which is non-radial, it is natural to expect the improvement of the Strichartz estimates by imposing radial assumption. In fact, such improvement is given by the work of J. Sterbenz \cite{sterbenz2} and Y. Cho - S. Lee \cite{choslee}. Hence one can use wider range of admissible Strichartz pairs. For instance, one may use
$$
\|e^{\mp it\langle\nabla\rangle}P_1f\|_{L^2_tL^{4+\epsilon}_x} \lesssim \|\langle\Omega\rangle^{\sigma}P_1f\|_{L^2_x}
$$
for arbitrarily small $\epsilon>0$ and $\sigma$ close to $\frac12$.
Nevertheless, we are not only aiming to improve the Sobolev index by the previous works but also interested in the low regularity problem with respect to the angular variables, namely, $0<\sigma\ll1$. For this purpose, we exploit the almost orthogonal decomposition by angular sectors together with cubes and then apply angular concentration estimates Lemma \ref{ang-con}. In this process, one may observe that the Low$\times$High frequency interactions become more difficult than the High$\times$High interaction, especially when the low frequency controls the angular frequency. Consequently we get a slightly bigger bound in this interactions. Fortunately, the Yukawa potential plays a distinguished role. Indeed, the potential is nothing but the Fourier multiplier with symbol $(b^2 + |\xi|^2)^{-1}$ and hence it is no harm to the summation.

In this manner we can prove Proposition \ref{bi-est} in the low-modulation regime. For high-modulation, we divide it into two cases.
\subsubsection*{High modulation - Low frequency}
As we have mentioned earlier, the situation when the modulation is larger than the highest-input frqeuncy is rather easier. We only consider when the modulation ranges from the lowest frequency to the highest frequency. As low-modulation regime, we deal with the High$\times$High and Low$\times$High frequency interactions. In the High$\times$High region, we can still use the null structure. However, we do not exploit the angular sector decomposition and hence angular regularity. Instead, one major observation is to decompose the modulation. In fact, we have the following decomposition:
\begin{align*}
\|C_dP_{\lambda_0}(C_{\le d}^{\pm_1}P_{\lambda_1}\varphi)^\dagger\beta(C^{\pm_2}_{\le d}P_{\lambda_2}\phi)\|_{L^2_tL^2_x} &\lesssim \|C_dP_{\lambda_0}(C^{\pm_1}_{\approx d}P_{\lambda_1}\varphi)^\dagger\beta(C^{\pm_2}_{\ll d}P_{\lambda_2}\phi)\|_{L^2_tL^2_x} \\
&\qquad\qquad +\|C_dP_{\lambda_0}(C^{\pm_1}_{\ll d}P_{\lambda_1}\varphi)^\dagger\beta(C^{\pm_2}_{\approx d}P_{\lambda_2}\phi)\|_{L^2_tL^2_x}.	
\end{align*}
Here $C_d^\pm$ is the modulation localisation operator. This can be easily derived by the support condition. The advantage of this observation is to allow the use of bound of high-modulation regime \eqref{bdd-high-mod}, which yields $d^{-\frac12}$ and this is truly helpful, since we are concerned with $d\gtrsim\min\{\lambda_0,\lambda_1,\lambda_2\}$. We will exploit the orthogonal decomposition by cubes as low-modulation case, and apply $L^4$-Strichartz estimates together with the bound \eqref{bdd-high-mod} to obtain the desired estimates. Even though we do not use the null structure in the Low$\times$High interaction, we can follow the aforementioned process and get the required bound.



\subsection*{Non-scattering for Coulomb potential}
We show a non-existence of scattering in $L^2_x$ for \eqref{main-eq} with the Coulomb potential $V_0(x) = \frac1{4\pi|x|}$.
In view of \cite{chagla, oz, chle}, there are trivial scattering conditions of \eqref{main-eq}. However, as observed in \cite{chle} for 2D problem, such scattering cannot occur for a certain class of solutions. 
To be precise, let us define $\mathcal I$ as follows:
$$
\mathcal I(\psi,\phi)(t):= \int_{\mathbb R^3} \left[V_0 *(\psi^\dagger  \phi)\right] (t,x) \left(\psi^\dagger \phi \right)(t,x)\, dx.
$$
Now we present our second theorem concerning non-scattering.
\begin{thm}\label{nonscatter-thm}
Assume that $\psi$ be a smooth solution to \eqref{main-eq} with $b=0$ which scatters in $L_x^2$ to a smooth solution $\psi_\infty^\ell$. If there exist $0<c<1$ and $t_*>0$ such that $\psi_\infty^\ell$ satisfies
\begin{align}\label{nonscatter-condi}
	\Big|\mathcal I(\psi_\infty^\ell,\beta\psi_\infty^\ell)(t) \Big| \ge c \mathcal I(\psi_\infty^\ell,\psi_\infty^\ell)(t)
\end{align}
for  $t>t_*$, then $\psi,\psi_\infty^\ell=0$ in $L_x^2$.
\end{thm}
If $\psi$ scatters in $L_x^2$ to $\psi_{\infty}^{\ell}$ satisfying \eqref{nonscatter-condi}, then by Lemmas \ref{infty-esti-1} and \ref{infty-weight} below one can find $0 < c' < 0$ and $t_{**} > 0$ such that
\begin{align*}
\mathcal I(\psi, \beta \psi )(t) \ge c' \mathcal I(\psi, \psi)(t)
\end{align*}
for any $t > t_{**}$.

It is essential to handle the lower bound of functional $H(t)  =  {\rm Im} \left< \psi(t), \psi_{\infty}^{\ell}(t) \right>_{L_x^2} $ by $\|\psi_{\infty}^{\ell}(0)\|_{L_x^2}$. We will show
$$
\left|\frac{d}{dt}H(t)\right| \ge c'|\mathcal I(\psi_{\infty}^{\ell}, \beta\psi_{\infty}^{\ell})(t)| + o(t^{-1}).
$$
Due to the matrix $\beta$ the value of $\mathcal I(\psi_{\infty}^{\ell}, \beta\psi_{\infty}^{\ell})$ may vanish. To avoid this we assume \eqref{nonscatter-condi} and hence obtain
$|\frac{d}{dt}H(t)| \gtrsim t^{-1}\|\psi_{\infty}^{\ell}(0)\|_{L_x^2}^2 + o(t^{-1})$ for sufficiently large $t$.
Therefore, if $\|\psi_{\infty}^{\ell}(0)\|_{L_x^2} > 0$, then the lower bound eventually will lead us to contradiction to the uniform boundedness of $H(t)$.

\subsection*{Organisation}
This paper is organized as follows. We introduce notations and preliminary setup in Section 2, which consists of Fourier localisation operators, basic analysis on the unit sphere, properties of $U^p-V^p$ spaces, and linear estimates related to the wave-type Strichartz estimates. In Section 3, we discuss the proof of Theorem \ref{gwp} and reduction to multilinear estimates. In Section 4, we present frequency-localised multilinear estimates (Proposition \ref{bi-est}), which play a crucial role in the proof of Theorem \ref{gwp}. 
Then Section 5 is devoted to the proof of Proposition \ref{bi-est}. Finally, we give the proof of non-scattering result in Section 6.

\subsection*{Notation}
As usual different positive constants, which are independent of dyadic numbers $\mu,\lambda$, and $d$ are denoted by the same letter $C$, if not specified. The inequalities $A \lesssim B$ and $A \gtrsim B$ means that $A \le CB$ and
$A \ge C^{-1}B$, respectively for some $C>0$. By the notation $A \approx B$ we mean that $A \lesssim B$ and $A \gtrsim B$, i.e., $\frac1CB \le A\le CB $ for some absolute constant $C$. We also use the notation $A\ll B$ if $A\le \frac1CB$ for some large constant $C$. Thus for quantites $A$ and $B$, we can consider three cases: $A\approx B$, $A\ll B$ and $A\gg B$. In fact, $A\lesssim B$ means that $A\approx B$ or $A\ll B$.

The spatial and space-time Fourier transform are defined by
$$
\widehat{f}(\xi) = \int_{\mathbb R^3} e^{-ix\cdot\xi}f(x)\,dx, \quad \widetilde{u}(\tau,\xi) = \int_{\mathbb R^{1+3}}e^{-i(t\tau+x\cdot\xi)}u(t,x)\,dtdx.
$$
We also write $\mathcal F_x(f)=\widehat{f}$ and $\mathcal F_{t, x}(u)=\widetilde{u}$. We denote the backward and forward wave propagation of a functiom $f$ on $\mathbb R^3$ by
$$
e^{\mp it\langle\nabla\rangle}f = \int_{\mathbb R^3}e^{ix\cdot\xi}e^{\mp it\langle\xi\rangle}\widehat{f}(\xi)\,d\xi.
$$
Finally, we shall use usual inner products for the normed space $\mathbb C^4$ and $L^2_x(\mathbb R^3)$. Namely, we write
$$
\langle\psi,\phi\rangle_{\mathbb C^4} = \psi^\dagger\phi,
$$
for spinor fields $\phi,\psi:\mathbb R^{1+3}\rightarrow\mathbb C^4$. We also write
$$
\langle f,g\rangle_{L^2_x} = \int_{\mathbb R^3} f(x)\overline{g(x)}\,dx,
$$
for any $L^2_x$-functions $f$ and $g$.
\section{Preliminary setup}
\subsection{Multipliers}\label{multi}
We fix a smooth function $\rho\in C^\infty_0(\mathbb R)$ such that $\rho$ is supported in the set $\{ \frac12<t<2\}$ and we let
$$
\sum_{\lambda\in2^{\mathbb Z}}\rho\left(\frac t\lambda\right) =1,
$$
and write $\rho_1=\sum_{\lambda\le1}\rho(\frac t\lambda)$ with $\rho_1(0)=1$. We define $\mathcal Q_\mu$ to be a finitely overlapping collection of cubes of diameter $\frac{\mu}{1000}$ covering $\mathbb R^3$, and let $\{ \rho_q\}_{q\in\mathcal Q_\mu}$ be a corresponding subordinate partition of unity. Now we define the standard Littlewood-Paley multipliers, for $\lambda\in 2^{\mathbb N}$, $\lambda>1$, $q\in\mathcal Q_\mu$, $d\in 2^{\mathbb Z}$:
$$
P_\lambda = \rho\left(\frac{|-i\nabla|}{\lambda}\right),\quad P_1=\rho_1(|-i\nabla|),\quad P_q = \rho_q(-i\nabla),\quad C^{\pm}_d = \rho\left(\frac{|-i\partial_t\pm\langle\nabla\rangle|}{d}\right).
$$
We also define $C^\pm_{\le d}=\sum_{\delta\le d}C^\pm_\delta$ and $C^\pm_{\ge d}$ is defined in the similar way. For simplicity we also write $C^+_d = C_d$.
Given $0<\alpha\lesssim1$, we define $\mathcal C_\alpha$ to be a collection of finitely overlapping caps of radius $\alpha$ on the sphere $\mathbb S^2$. If $\kappa\in\mathcal C_\alpha$, we let $\omega_\kappa$ be the centre of the cap $\kappa$. Then we define $\{\rho_\kappa\}_{\kappa\in\mathcal C_\alpha}$ to be a smooth partition of unity subordinate to the conic sectors $\{ \xi\neq0 , \frac{\xi}{|\xi|}\in\kappa \}$ and denote the angular Fourier localisation multipliers by
$
R_\kappa = \rho_\kappa(-i\nabla).
$
\subsection{Analysis on the sphere}\label{an-sph}
We introduce some basic facts from harmonic analysis on the unit sphere. The most of ingredients can be found in \cite{candyherr,sterbenz1}. We also refer the readers to \cite{steinweiss} for more systematic introduction to the spherical harmonics. We let $Y_{\ell}$ be the set of homogeneous harmonic polynomial of degree $\ell$. Then define $\{ y_{\ell,n} \}_{n=0}^{2\ell}$ a set of orthonormal basis for $Y_{\ell}$, with respect to the inner product:
\begin{align}
\langle y_{\ell,n},y_{\ell',n'}\rangle_{L^2_\omega(\mathbb S^2)} = \int_{\mathbb S^2}{y_{\ell,n}(\omega)} \overline{y_{\ell',n'}(\omega)}\,d\omega.
\end{align}
Given $f\in L^2_x(\mathbb R^3)$, we have the orthogonal decomposition as follow:
\begin{align}
f(x) = \sum_{\ell}\sum_{n=0}^{2\ell}\langle f(|x|\omega),y_{\ell,n}(\omega)\rangle_{L^2_\omega(\mathbb S^2)}y_{\ell,n}\big(\frac{x}{|x|}\big).
\end{align}
For a dyadic number $N>1$, we define the spherical Littlewood-Paley decompositions by
\begin{align}
H_N(f)(x) & = 	\sum_{\ell}\sum_{n=0}^{2\ell}\rho\left(\frac\ell N\right)\langle f(|x|\omega),y_{\ell,n}(\omega)\rangle_{L^2_\omega(\mathbb S^2)}y_{\ell,n}\big(\frac{x}{|x|}\big), \\
H_1(f)(x) & = \sum_{\ell}\sum_{n=0}^{2\ell}\rho_{\le1}(\ell)\langle f(|x|\omega),y_{\ell,n}(\omega)\rangle_{L^2_\omega(\mathbb S^2)}y_{\ell,n}\big(\frac{x}{|x|}\big).
\end{align}
Since $-\Delta_{\mathbb S^2}y_{\ell, n} = \ell(\ell+1)y_{\ell, n}$, by orthogonality one can readily get
$$\|\langle\Omega\rangle^\sigma f\|_{L^2_\omega({\mathbb S^2})} \approx \left\|\sum_{N\in2^{\mathbb N}\cup\{0\}}N^\sigma H_Nf\right\|_{L^2_\omega({\mathbb S^2})}.$$


\begin{lem}[Lemma 7.1. of \cite{candyherr}]\label{sph-ortho}
Let $N\ge1$. Then $H_N$ is uniformly bounded on $L^p(\mathbb R^3)$ in $N$, and $H_N$ commutes with all radial Fourier multipliers. Moreover, if $N'\ge1$, then either $N\approx N'$ or
$$
H_N\Pi_\pm H_{N'}=0.
$$	
\end{lem}
By Lemma \ref{sph-ortho}, we see that $H_N$ commutes with the $P_\lambda$ and $C^\pm_d$ multipliers since we can write $C_d^\pm=e^{\mp it\langle\nabla\rangle}\rho(-\frac{i\partial_t}{d})e^{\pm it\langle\nabla\rangle}$. On the other hand, we note that $H_N$ does not commute with the cube and cap localisation operators $R_\kappa$ and $P_q$, which are non-radial.

\subsection{Adapted function spaces}\label{ftn-sp}

Let $1\le p\le+\infty$ and $\mathcal I=\left\{ \{t_k\}_{k=0}^K : t_k\in\mathbb R, t_k<t_{k+1} \right\}$ be the set of increasing sequences of real numbers.
We define the $p$-variation of $v$ to be
$$
|v|_{V^p} = \sup_{ \{t_k\}_{k=0}^K\in\mathcal I } \left( \sum_{k=0}^K\|v(t_k)-v(t_{k-1})\|_{L^2_x}^p \right)^\frac1p
$$
Then the Banach space $V^p$ can be defined to be all right continuous functions $v:\mathbb R\rightarrow L^2_x$ such that the quantity
$$
\|v\|_{V^p} = \|v\|_{L^\infty_tL^2_x} + |v|_{V^p}
$$
is finite. Set $\|u\|_{V^p_\pm}=\|e^{\mp it\langle\nabla\rangle}u\|_{V^p}$. We recall basic properties of $V^2_\pm$ space from \cite{candyherr, candyherr1}. For more about $U^p-V^p$ space see \cite{haheko}.

The following lemma is on a simple bound in the high-modulation region.
\begin{lem}[Corollary 2.18. of \cite{haheko}]
Let $2\le q\le\infty$. For $d\in2^{\mathbb Z}$, we have
\begin{align}\label{bdd-high-mod}
\begin{aligned}
\|C^{\pm}_du\|_{L^q_tL^2_x} \lesssim d^{-\frac1q}\|u\|_{V^2_\pm},\\
\|C^{\pm}_{\ge d}u\|_{L^q_tL^2_x} \lesssim d^{-\frac1q}\|u\|_{V^2_\pm}.
\end{aligned}
\end{align}
\end{lem}
Now we present an energy inequality. See also \cite[Proposition 2.10]{haheko}.
\begin{lem}[Lemma 7.3. of \cite{candyherr}]\label{lem-v-dual}
Let $F\in L^\infty_tL^2_x$, and suppose that
$$
\sup_{\|P_\lambda H_Nv\|_{V^2_\pm}\lesssim1}\left|\int_{\mathbb R} \langle P_\lambda H_Nv(t),F(t)\rangle_{L^2_x}\,dt \right| <\infty.
$$	
If $u\in C(\mathbb R,L^2_x)$ satisfies $-i\partial_tu\pm\langle\nabla\rangle u=F$, then $P_\lambda H_Nu\in V^2_\pm$ and we have the bound
\begin{align}\label{v-dual}
\|P_\lambda H_Nu\|_{V^2_\pm} \lesssim \|P_\lambda H_Nu(0)\|_{L^2_x} + \sup_{\|P_\lambda H_Nv\|_{V^2_\pm}\lesssim1}\left|\int_{\mathbb R} \langle P_\lambda H_Nv(t),F(t)\rangle_{L^2_x}\,dt \right|.
\end{align}
\end{lem}

We recall the uniform disposability of the modulation cutoff multipliers, which reads for $1\le q,r\le\infty$,
\begin{align}\label{uni-dis1}
\|C^\pm_{\le d}P_\lambda R_\kappa u\|_{L^q_tL^r_x} + \|C^\pm_dP_{\lambda}R_\kappa u\|_{L^q_tL^r_x} \lesssim \|P_\lambda R_\kappa u\|_{L^q_tL^r_x},
\end{align}
if $\kappa\in\mathcal C_\alpha,\ d\gtrsim\alpha^2\lambda,$ and $\alpha\gtrsim\lambda^{-1}$.
Since convolution with $L^1_t(\mathbb R)$ functions is bounded on the $V^2$ space, we also have for every $d\in2^{\mathbb Z}$,
\begin{align}\label{uni-dis}
\|C^\pm_{\le d}u\|_{V^2_\pm} \lesssim \|u\|_{V^2_\pm}.
\end{align}
To prove the scattering result, we shall use the following lemma:
\begin{lem}[Lemma 7.4. of \cite{candyherr}]\label{v-scatter}
Let $u\in V^2_\pm$. Then there exists $f\in L^2_x$ such that $\|u(t)-e^{\mp it\langle\nabla\rangle}f\|_{L^2_x}\rightarrow0$ as $t\rightarrow\pm\infty$.
\end{lem}

\subsection{Auxiliary estimates}
In this section we provide the key ingredients to the proof of Theorem \ref{gwp}. We introduce the null-form-bound, localised Strichartz estimates, and various square-summation. We would like to highlight that there is nothing new, and hence we simply list several estimates used in the sequel without proof. However, we also encourage readers to read the reference \cite{bejeherr,candyherr,sterbenz1,sterbenz2}.

It is well-known fact that the nonlinearity in the system \eqref{dirac-eq} has null structure, which gives the cancellation property. To reveal null form, we write
\begin{align}\label{dec-di}
\begin{aligned}
(\Pi_{\pm_1}\phi)^\dagger\beta\Pi_{\pm_2}\varphi = &[(\Pi_{\pm_1}-\Pi_{\pm_1}(x))\phi]^\dagger\beta\Pi_{\pm_2}\varphi+(\Pi_{\pm_1}\phi)^\dagger\beta(\Pi_{\pm_2}-\Pi_{\pm_2}(y))\varphi \\
&\qquad\qquad+\phi^\dagger\Pi_{\pm_1}(x)\beta\Pi_{\pm_2}(y)\varphi,
\end{aligned}
\end{align}
for any $x,y\in\mathbb R^3$. Then we have the following null-form-type bound:
\begin{align}\label{dirac-null}
|\Pi_{\pm_1}(\xi)\beta\Pi_{\pm_2}(\eta)| \lesssim \angle(\pm_1\xi,\pm_2\eta)+\frac{|\pm_1|\xi|\pm_2|\eta||}{\langle\xi\rangle\langle\eta\rangle}.
\end{align}
To exploit the null form for the first and second terms of \eqref{dec-di}, we use the following lemma:
\begin{lem}[Lemma 8.1. of \cite{candyherr}]\label{lem-null}
Let $1<r<\infty$. If $\lambda\ge1$, $\alpha\gtrsim\lambda^{-1},\ \kappa\in\mathcal C_\alpha$, then
$$
\|(\Pi_\pm-\Pi_{\pm}(\lambda\omega(\kappa)))R_\kappa P_\lambda f\|_{L^r_x} \lesssim \alpha\|R_\kappa P_\lambda f\|_{L^r_x}.
$$	
\end{lem}
We recall Strichartz estimates for wave equation, which plays a significant role in the proof of Theorem \ref{gwp}. Note that an additional angular regularity allows to extend admissible Strichartz pairs.
\begin{lem}[Lemma 3. of \cite{candyherr1}]\label{wave-stri}
Let $2<q\le\infty$. If $0<\mu\le\lambda,\ N\ge1,$ and $\frac1q+\frac1r=\frac12$, then for every $q\in\mathcal Q_\mu$ we have
$$
\|e^{\mp it\langle\nabla\rangle}P_qP_\lambda f\|_{L^q_tL^r_x} \lesssim \mu^{\frac12-\frac1r}\lambda^{\frac12-\frac1r}\|P_qP_\lambda f\|_{L^2_x}.
$$	
Moreover, by spending additional angular regularity, if $\frac1q+\frac2r<1$, then we have for $\epsilon>0$
$$
\|e^{\mp it\langle\nabla\rangle}P_\lambda H_Nf\|_{L^q_tL^r_x} \lesssim \lambda^{3(\frac12-\frac1r)-\frac1q}N^{\frac12+\epsilon}\|P_\lambda H_Nf\|_{L^2_x}.
$$
\end{lem}
For the first estimate, see Lemma 3.2. of \cite{cyang}, or Lemma 3.1, of \cite{bejeherr}. The second estimate can be found in Theorem 1. of \cite{choslee}. Note that the use of decomposition by smaller cubes $q\in\mathcal Q_\mu$ gives us better estimates. However, we need to recover the term $\|P_\lambda f\|_{L^2_x}$ by square-summation with respect to cubes $q\in\mathcal Q_\mu$ and this would result in a certain loss in estimates. Fortunately, the following lemma allows this loss to be absorbed elsewhere.
\begin{lem}[Lemma 8.6. of \cite{candyherr}]\label{square-sum}
Let $\{P_j\}_{j\in\mathcal J}$ and $\{M_j\}_{j\in\mathcal J}$ be a collection of spatial Fourier multipliers. Suppose that the symbols of $P_j$ have finite overlap, and
$$
\|M_jP_jf\|_{L^2_x} \lesssim \delta \|P_jf\|_{L^2_x}
$$	
for some $\delta>0$. Let $q>2,\ r\ge2$. Suppose that there exists $A>0$ such that for every $j$ we have the bound
$$
\|e^{\mp it\langle\nabla\rangle}P_jf\|_{L^q_tL^r_x} \le A\|P_jf\|_{L^2_x}.
$$
Then for every $\epsilon>0$, we have
$$
\left(\sum_{j\in\mathcal J}\|M_jP_jv\|_{L^q_tL^r_x}^2\right)^\frac12 \lesssim \delta|\mathcal J|^\epsilon A\|v\|_{V^2_\pm}.
$$
Here $|\mathcal J|$ is the cardinal number of the set $\mathcal J$.
\end{lem}
For future use, we list some of the direct results of Lemma \ref{square-sum} as follows:
\begin{lem}
Let $1\le\mu\lesssim\lambda,\ \alpha\gtrsim\lambda^{-1},\ \epsilon>0$. For $q,r$ satisfying the condition as in Lemma \ref{wave-stri}, we have
\begin{align}
\left( \sum_{q\in\mathcal Q_\mu}\sum_{\kappa\in\mathcal C_\alpha}\|R_\kappa P_q u_{\lambda,N}\|^2_{L^q_tL^r_x} \right)^\frac12 & \lesssim \alpha^{-\epsilon}\left(\frac\mu\lambda\right)^{-\epsilon}(\mu\lambda)^{\frac12-\frac1r}\|u_{\lambda,N}\|_{V^2_\pm},\label{cube-lin} \\
\left(\sum_{\kappa\in\mathcal C_\alpha}\|R_\kappa u_{\lambda,N}\|_{L^q_tL^r_x}^2 \right)^\frac12 & \lesssim \alpha^{-\epsilon}\lambda^{3(\frac12-\frac1r)-\frac1q}N^{\frac12+\epsilon}\|u_{\lambda,N}\|_{V^2_\pm}. \label{ang-lin}
\end{align}
Here we write $P_\lambda H_N u=u_{\lambda,N}$.
\end{lem}
So far we have seen that the additional angular regularity gives rise to the improvement of space-time Strichartz estimates. However, we are interested in only small amount of angular regularity, namely, $\sigma\ll1$. To achieve this low-regularity-condition, we exploit the following so-called angular concentration estimates:
\begin{lem}[Lemma 8.5. of \cite{candyherr}]\label{ang-con}
Let $2\le p<\infty$, and $0\le s<\frac2p$. If $\lambda,N\ge1$, $\alpha\gtrsim\lambda^{-1}$, and $\kappa\in\mathcal C_\alpha$, then we have
$$
\|R_\kappa P_\lambda H_N f\|_{L^p_x(\mathbb R^3)} \lesssim (\alpha N)^s \|P_\lambda H_N f\|_{L^p_x(\mathbb R^3)}.
$$	
\end{lem}
The proof can be found in Lemma 5.2. of \cite{sterbenz2}.

\subsection{General resonance identity}
In the proof of Theorem \ref{gwp}, we shall follow the standard fixed-point argument. That is, by the use of Lemma \ref{lem-v-dual}, we consider the following quartilinear form:
$$
\int_{\mathbb R^{1+3}}(\varphi^\dagger\beta\phi)(\uppsi^\dagger\beta\psi)\,dtdx.
$$
After using the H\"{o}lder inequality, we need to deal with bilinear estimates, which has essentially trilinear expression via $L^2$-duality:
\begin{align}\label{tri-int}
\int_{\mathbb R^{1+3}}\phi\psi^\dagger\beta\varphi\, dxdt.
\end{align}
Suppose that $\phi,\psi,\varphi$ have small modulation. We also assume that the Fourier transform of $\phi$ is supported in $\{ |\tau+\langle\xi\rangle|\le d \}$, the support of $\widetilde{\psi}$ is contained in the set $\{ |\tau\pm_1\langle\xi\rangle|\le d \}$, and $\widetilde{\varphi}$ is supported in the set $\{ |\tau\pm_2\langle\xi\rangle|\le d \}$ for some $d\in 2^{\mathbb Z}$. Then the integral \eqref{tri-int} vanishes unless
$$
|\langle\xi-\eta\rangle\mp_1\langle\xi\rangle\pm_2\langle\eta\rangle|\lesssim d.
$$
Now we define the modulation function:
$$
\mathcal M_{\pm_1,\pm_2}(\xi,\eta) = |\langle\xi-\eta\rangle\mp_1\langle\xi\rangle\pm_2\langle\eta\rangle|.
$$
We first note the symmetry properties of $\mathcal M_{\pm_1,\pm_2}$, for example, we have $\mathcal M_{+,+}(\xi,\eta)=\mathcal M_{-,-}(\eta,\xi)$, and $\mathcal M_{\pm,\mp}(\xi,\eta)=\mathcal M_{\pm,\mp}(\eta,\xi)$.
\begin{lem}[Lemma 8.7. of \cite{candyherr}]\label{resonance}
We have
\begin{align*}
\mathcal M_{-,+}(\xi,\eta) &\gtrsim \langle\xi\rangle+\langle\eta\rangle, \\
\mathcal M_{\pm,\pm}(\xi,\eta) &\gtrsim \frac{1}{\langle\xi-\eta\rangle}\left( \frac{(|\xi|-|\eta|)^2}{\langle\xi\rangle\langle\eta\rangle}+|\xi||\eta|\angle(\xi,\eta)^2+1 \right), \\
\mathcal M_{-,-}(\xi,\eta) &\gtrsim \frac{|\xi-\eta||\xi|}{\langle\xi\rangle+\langle\eta\rangle}\angle(\xi-\eta,-\xi)^2, \\
\mathcal M_{+,+}(\xi,\eta) &\gtrsim \frac{|\xi-\eta||\eta|}{\langle\xi\rangle+\langle\eta\rangle}\angle(\xi-\eta,\eta)^2.
\end{align*}
\end{lem}

\section{Proof of Theorem \ref{gwp}}

In this section we prove Theorem \ref{gwp} via standard fixed-point argument. To be precise, we shall find the solution to \eqref{main-eq} in a complete metric space $(F^\sigma(\delta), d)$ defined as
\begin{align*}
F^{\sigma}(\delta) &:= \Big\{ \psi : \|\psi\|_{F^{\sigma}}:= \|\psi_+\|_{F_+^{\sigma}} + \|\psi_-\|_{F_-^{\sigma}} < \delta  \Big\},\quad d(\psi, \phi) := \|\psi - \phi\|_{F^\sigma},
\end{align*}
where $$\|u\|_{F_\pm^{ \sigma}} : = \left(\sum_{\lambda, N \in 2^{\mathbb N\cup \{0\}}} N^{2\sigma}\|P_\lambda H_N u\|_{V_\pm^2}^2\right)^\frac12.$$
We also define the map $\mathcal Y$ on $F^\sigma(\delta)$ by
\begin{align}\label{map}
\mathcal Y(\psi) := \sum_{\pm_0 \in \{ \pm \}} e^{-\pm_0 i t \langle \nabla \rangle} \Pi_{\pm_0}\psi_0  +i \sum_{\pm_j;j=0,1,2,3}   \mathcal N_{\pm_0} (\psi_{1},\psi_{2},\psi_{3})(t),
\end{align}
where
\begin{align*}
	\mathcal N_{\pm_0}(\psi_{1},\psi_{2},\psi_{3})(t) = \int_0^t e^{\mp_0 i(t-t')\langle \nabla \rangle}\Pi_{\pm_0}[(V*\psi_1^\dagger\beta\psi_2)\beta\psi_{3}](t')\, dt'.
\end{align*}
Here we put $\psi_{j} := \Pi_{\pm_j}\psi$. Then we need to show that $\mathcal Y$ is a contraction mapping on $F^\sigma(\delta)$. Indeed the linear part of \eqref{map} can be handled as follows:
\begin{align}\label{linear-est}
	\left\|e^{\mp_0 it\langle \nabla \rangle} \Pi_{\pm_0} \psi_0 \right\|_{F_{\pm_0}^{\sigma}}^2 = \sum_{\lambda \ge 1}\sum_{N \ge 1} N^{2\sigma}\left\|\chi_{[0,\infty)} P_\lambda H_N \Pi_{\pm_0} \psi_0 \right\|_{V_{\pm_0}^2}^2 \le C\|\langle \Omega \rangle^\sigma \psi_0\|_{L_x^2}^2
\end{align}
We are left to control the nonlinearity $\mathcal N_{\pm_0}(\psi_1,\psi_2,\psi_3)$ of \eqref{map}.

\begin{prop}\label{main-est}
	Let $\sigma>0$ and $\psi_j \in F_{\pm_j}^{\sigma}(j=1,2,3)$. Then we have the following multilinear estimates:
	\begin{align}
		\Big\| \mathcal N_{\pm_0} (\psi_1,\psi_2,\psi_3) \Big\|_{F^\sigma_{\pm_0}} & \lesssim  \prod_{j=1}^3 \|\psi_j\|_{F^\sigma_{\pm_j}}.
	\end{align}
\end{prop}
Once Proposition \ref{main-est} has been proved, this trilinear estimate together with linear estimate \eqref{linear-est} will lead us that
\begin{align*}
	\|\mathcal Y(\psi)\|_{F^{\sigma}} = \sum_{\pm_0}\|\Pi_{\pm_0}\mathcal Y(\psi)\|_{F_{\pm_0}^{\sigma}} \le C \left( \|\langle \Omega \rangle^\sigma\psi_0\|_{L_x^2} + \|\psi\|_{F^\sigma}^3 \right).
\end{align*}
If $\delta$ is small enough so that $C\delta^3 \le \frac\delta8$ and the initial data $\psi_0$ satisfies the smallness condition, namely, $C \left\|\langle \Omega \rangle^\sigma\psi_0 \right\|_{L_x^2} \le \frac\delta2$,
then $\mathcal Y$ is a self-mapping on $F^\sigma(\delta)$. Furthermore, we get
\begin{align*}
	\|\mathcal{Y}(\psi) - \mathcal{Y}(\phi)\|_{F^\sigma} &\le C\left(\|\psi\|_{F^\sigma}+ \|\phi\|_{F^\sigma}\right)^2 \|\psi- \phi\|_{F^\sigma} \le 4C\delta^2\|\psi- \phi\|_{F^\sigma} \le \frac12\|\psi- \phi\|_{F^\sigma}.
\end{align*}
Hence $\mathcal Y: F^\sigma(\delta) \to F^\sigma(\delta)$ is a contraction mapping for sufficiently small $\delta$, which completes the proof of global well-posedness of \eqref{main-eq}.

Now we move onto the scattering property of \eqref{main-eq}. We have the solutions $\psi_\pm\in F^\sigma_\pm$, and hence $\langle\Omega\rangle^\sigma\psi_\pm\in V^2_\pm$. By Lemma \ref{v-scatter}, there exists function $\varphi_\pm \in L^{2, \sigma}_x:=\langle\Omega\rangle^{-\sigma}L^2_x$ such that
$$
\|\psi_\pm - e^{\mp it\langle\nabla\rangle}\varphi\|_{L^{2, \sigma}_x}\rightarrow 0,
$$
as $t\rightarrow\pm\infty$, which completes the proof of scattering of \eqref{main-eq}.
Finally, we are left to show the multilinear estimates Proposition \ref{main-est}.



\section{Proof of Proposition \ref{main-est}}

In the remainder of this paper, we focus on the proof of multilinear estimates Proposition \ref{main-est}. First, after an application of Lemma \ref{lem-v-dual}, we write
\begin{align*}
	&\Big\|\mathcal N_{\pm_4}(\varphi,\phi,\psi) \Big\|_{V^2_{\pm_4}}\\
	& \lesssim \sup_{\|P_{\lambda_4} H_{N_4}\uppsi\|_{V^2_{\pm_4}}\lesssim1}\left| \int_{\mathbb R^{1+3}} (P_{\lambda_4} H_{N_4}\uppsi)^\dagger	\Pi_{\pm_4} (V*(\varphi^\dagger\beta\phi))\beta\psi\,dtdx\right| \\
	& \lesssim \sup_{\|P_{\lambda_4} H_{N_4}\uppsi\|_{V^2_{\pm_4}}\lesssim1}\left|\int_{\mathbb R^{1+3}} V*(\varphi^\dagger\beta\phi)(\Pi_{\pm_4} \uppsi_{\lambda_4,N_4})^\dagger\beta\psi \,dtdx\right| \\
	& \lesssim \sup_{\|\uppsi_{\lambda_4,N_4}\|_{V^2_{\pm_4}}\lesssim1}\left( \sum_{\lambda_1,\lambda_2,\lambda_3\ge1}\sum_{N_1,N_2,N_3\ge1}\left|\int_{\mathbb R^{1+3}} V*(\varphi^\dagger_{\lambda_1,N_1}\beta\phi_{\lambda_2,N_2})(\uppsi^\dagger_{\lambda_4,N_4}\beta\psi_{\lambda_3,N_3})\,dtdx \right| \right),
\end{align*}
where we put $\psi_{\lambda,N}=\Pi_\pm P_\lambda H_N\psi$. From now on, we turn our attention to the following quardrilinear form:
\begin{align}\label{main-int}
	\mathcal I_{\lambda,N} := \left|\int_{\mathbb R^{1+3}}V*(\varphi^\dagger_{\lambda_1,N_1}\beta\phi_{\lambda_2,N_2})(\uppsi^\dagger_{\lambda_4,N_4}\beta\psi_{\lambda_3,N_3})\,dtdx\right|.
\end{align}
We first note that if $\xi_j$ are the spatial Fourier frequencies for the functions in the integrand, by Plancherel's theorem the integral \eqref{main-int} vanishes unless
$$
-\xi_1+\xi_2+\xi_3-\xi_4=0.
$$
If $\xi_0$ is the output frequency of the bilinear form $\varphi^\dagger\beta\phi$, then we have
$$
\xi_0 = -\xi_1+\xi_2 = \xi_4-\xi_3.
$$
Thus, if the output frequency $\xi_0$ is localised in an annulus of dyadic radius $\lambda_0$, the standard Littlewood-Paley trichotomy must give the following frequency interactions:
\begin{align}\label{fre-re}
	\begin{aligned}
		\min\{\lambda_0,\lambda_1,\lambda_2\} \lesssim \textrm{med}\{\lambda_0,\lambda_1,\lambda_2\} \approx \max\{\lambda_0,\lambda_1,\lambda_2\}, \\
		\min\{\lambda_0,\lambda_3,\lambda_4\} \lesssim \textrm{med}\{\lambda_0,\lambda_3,\lambda_4\} \approx \max\{\lambda_0,\lambda_3,\lambda_4\}.
	\end{aligned}
\end{align}
In view of Lemma \ref{sph-ortho}, we also have the similar relation for the angular frequencies $N_j$, $j=0,1,\cdots,4$:
\begin{align}\label{angfre-re}
	\begin{aligned}
		\min\{ N_0,N_1,N_2 \} & \lesssim \textrm{med}\{ N_0,N_1,N_2 \} \approx \max\{ N_0,N_1,N_2 \}, \\
		\min\{ N_0,N_3,N_4 \} & \lesssim \textrm{med}\{ N_0,N_3,N_4 \} \approx \max\{ N_0,N_3,N_4 \}.
	\end{aligned}	
\end{align}
We also consider the relation between spatial frequency $\lambda_j$ and modulation $d$ (distance to the cone). For this purpose we are concerned with three cases as follows:
\begin{align}
d \lesssim & \,\lambda_{\min}, \label{high-fre-low-mod} \\
\lambda_{\min} \ll d &\ll \lambda_{\max}, \label{low-fre-high-mod1} \\
\lambda_{\max}& \lesssim d. \label{low-fre-high-mod2}
\end{align}
The third case \eqref{low-fre-high-mod2} is rather easier than other cases.
Indeed, it suffices to consider one of spinor fields has higher modulation $\gtrsim\lambda_{\max}$. We may assume that $\varphi_{\lambda_1,N_1}$ has the highest modulation and write $\varphi_{\lambda_1,N_1} = C^{\pm_1}_{\ge d}\varphi_{\lambda_1,N_1}$. 
Then we consider the following quartilinear expression.
\begin{align*}
\mathcal I_{\lambda,N} & \lesssim \left|\int_{\mathbb R^{1+3}}  V*\Big( \left[ C_{\ge d}^{\pm_1}\varphi_{\lambda_1,N_1}\right]^\dagger\beta\phi_{\lambda_2,N_2}\Big)(\uppsi^\dagger_{\lambda_4,N_4}\beta\psi_{\lambda_3,N_3})\,dtdx\right|.	
\end{align*}
We first use H\"older's inequality.
\begin{align*}
\mathcal I_{\lambda,N} & \lesssim 	\langle {\lambda_0} \rangle^{-2}    \left\| P_{{\lambda_0}} \left[  \left(C_{\ge d}^{\pm_1}\varphi_{\lambda_1,N_1}\right)^\dagger\beta\phi_{\lambda_2,N_2}\right]\right\|_{L_t^2L_x^{1+\epsilon}}\left\|P_{{\lambda_0}} \left(\uppsi^\dagger_{\lambda_4,N_4}\beta\psi_{\lambda_3,N_3}\right)\right\|_{L_t^2L_x^{\frac{1+\epsilon}{\epsilon}}}.
\end{align*}
Then by the use of H\"older's inequality for $L^2_tL^{1+\epsilon}_x$ norm and Bernstein's inequality for $L^2_tL^{\frac{1+\epsilon}{\epsilon}}_x$, we see that
\begin{align*}
\mathcal I_{\lambda,N} & \lesssim \lambda_0^{-2}\|C^{\pm_1}_{\ge d}\varphi_{\lambda_1,N_1}\|_{L^2_{t,x}}\|\phi_{\lambda_2,N_2}\|_{L^\infty_tL^{\frac{2(1+\epsilon)}{1-\epsilon}}_x}\lambda_0^{\frac32\frac{1-\epsilon}{1+\epsilon}}\|P_{\lambda_0}(\uppsi^\dagger_{\lambda_4,N_4}\beta \psi_{\lambda_3,N_3})\|_{L^2_{t,x}}.	
\end{align*}
We use boundedness in high-modulation regime \eqref{bdd-high-mod} for $C^{\pm_1}_{\ge d}\varphi_{\lambda_1,N_1}$ to gain $d^{-\frac12}$. We also use Bernstein's inequality for $\phi_{\lambda_2,N_2}$ and then apply usual energy estimate. For the $L^2_{t,x}$-bilinear estimates, we refer to Corollary 3.8 of \cite{cyang}, which yields
$$
\|P_{\lambda_0}(\uppsi_{\lambda_4}^\dagger\beta\,\psi_{\lambda_3})\|_{L^2_{t,x}} \lesssim \lambda_0\|\psi_{\lambda_3}\|_{V^2_{\pm_3}}\|\uppsi_{\lambda_4}\|_{V^2_{\pm_4}}.
$$
 Hence we obtain
\begin{align*}
\sum_{d\gtrsim\lambda_{\max}}\mathcal I_{\lambda,N}& \lesssim	 \left(\frac{\lambda_0}{\lambda_{\max}}\right)^{\frac12-\frac{3\epsilon}{1+\epsilon}}\|\varphi_{\lambda_1,N_1}\|_{V^2_{\pm_1}}\|\phi_{\lambda_2,N_2}\|_{V^2_{\pm_2}}\|\psi_{\lambda_3,N_3}\|_{V^2_{\pm_3}}\|\uppsi_{\lambda_4,N_4}\|_{V^2_{\pm_4}}.
\end{align*}
Then usual square summation with respect to $\lambda_j$ and $N_j$ gives the required estimate as Proposition \ref{main-est}.

Consequently, we are left to consider the cases \eqref{high-fre-low-mod}, \eqref{low-fre-high-mod1}.
We use H\"{o}lder's inequality for $\mathcal I_{\lambda,N}$ to obtain
\begin{align*}
	\mathcal I_{\lambda,N} \lesssim \sum_{\lambda_0,N_0\ge1}\|\langle \nabla\rangle^{-2}P_{\lambda_0}H_{N_0}(\varphi^\dagger_{\lambda_1,N_1}\beta\phi_{\lambda_2,N_2})\|_{L^2_{t,x}}\|P_{\lambda_0}H_{N_0}(\uppsi^\dagger_{\lambda_4,N_4}\beta\psi_{\lambda_3,N_3})\|_{L^2_{t,x}},
\end{align*}
where $\langle\nabla\rangle$ is the Fourier multiplier, whose symbol is given by $\langle\xi\rangle$.
In the region $d\ll\lambda_{\max}$, we have the following frequency-localised $L^2$-bilinear estimates.
\begin{prop}\label{bi-est}
	Let $\epsilon>0$. For some $\delta>0$, we have the following estimates:
	\begin{align}
		\|P_{\lambda_0}H_{N_0}(\varphi_{\lambda_1,N_1}^\dagger\beta\phi_{\lambda_2,N_2})\|_{L^2_{t,x}} & \lesssim \lambda_0\left(\frac{\lambda_{\min}}{\lambda_{\max}}\right)^\delta (N_{\min})^\epsilon \|\varphi_{\lambda_1,N_1}\|_{V^2_{\pm_1}}\|\phi_{\lambda_2,N_2}\|_{V^2_{\pm_2}},
	\end{align}
	where $\lambda_{\min}$ and $\lambda_{\max}$ are the minimum and maximum of $\{\lambda_0,\lambda_1,\lambda_2\}$, respectively and $N_{\min}=\min\{N_0,N_1,N_2\}$.
\end{prop}

\subsection{Proof of Proposition \ref{main-est}}

We first prove Theorem \ref{main-est} restricted to the region $d\ll\lambda_{\max}$. Then we write
\begin{align*}
& \Big\|  \mathcal N_{\pm_4} (\varphi,\phi,\psi)\Big\|_{F^\sigma_{\pm_4}}^2 \\
& = \sum_{\lambda_4\ge1}\sum_{N_4\ge1} (N_4)^{2\sigma} \Big\|P_{\lambda_4}H_{N_4}	\mathfrak I^{\pm_4}[\Pi_{\pm_4}(V*(\varphi^\dagger\beta\phi))\beta\psi]\Big\|_{V^2_{\pm_4}}^2 \\
& \lesssim \sum_{\lambda_4,N_4}N_4^{2\sigma} \sup_{\|\uppsi_{\lambda_4,N_4}\|_{V^2_{\pm_4}}\lesssim1}\left( \sum_{\lambda_j,N_j}\left|\int_{\mathbb R^{1+3}} V*(\varphi^\dagger_{\lambda_1,N_1}\beta\phi_{\lambda_2,N_2})(\uppsi^\dagger_{\lambda_4,N_4}\beta\psi_{\lambda_3,N_3})\,dtdx \right| \right)^2 \\
 	 & \lesssim \sum_{\lambda_4,N_4}N_4^{2\sigma}\left( \sum_{\lambda_j,N_j}\left(\frac{\lambda_{\min}\lambda_{\min}'}{\lambda_{\max}^{12}\lambda_{\max}^{34}}\right)^\delta (N_{\min}N_{\min}')^\epsilon \|\varphi_{\lambda_1,N_1}\|_{V^2_{\pm_1}}\|\phi_{\lambda_2,N_2}\|_{V^2_{\pm_2}}\|\psi_{\lambda_3,N_3}\|_{V^2_{\pm_3}} \right)^2 \\
 	 & \quad =: I(\varphi,\phi, \psi),
 \end{align*}
 where we used H\"older's inequality and Proposition \ref{bi-est}.
 Here $\lambda_{\min},\lambda_{\min}'$ are the first and second lowest terms of $\lambda_1,\lambda_2,\lambda_3,\lambda_4$ and $N_{\min},N_{\min}'$ are defined in the similar way. Note that the factor $\left(\frac{\lambda_{\min}\lambda_{\min}'}{\lambda_{\max}^{12}\lambda_{\max}^{34}}\right)^\delta$ plays a role as kernel in the square summation with respect to the $\lambda_j$. By Young's convolution inequality for $\ell^p$ spaces we obtain the desired estimates as Proposition \ref{main-est}. We omit the details. See also Remark 1 of \cite{tes1}.

\section{Bilinear estimates: Proof of Proposition \ref{bi-est}}
In this section, we prove various frequency-localised bilinear estimates, which implies Proposition \ref{bi-est} when the modulation is bounded above by the highest frequency.

\subsection{High frequency - Low modulation}
First, we are concerned with low-modulation regime. The following bilinear estimates imply Proposition \ref{bi-est} when the modulation is less than the lowest frequency.
\begin{thm}
Let $\sigma>0$. For arbitrarily small $\epsilon>0$, we have
\begin{align}\label{est-hh1}
\|P_{\lambda_0}H_{N_0}(\varphi_{\lambda_1,N_1}^\dagger\beta\,\phi_{\lambda_2,N_2})	\|_{L^2_{t,x}} & \lesssim \lambda_0\left(\frac{\lambda_0}{\lambda_1}\right)^{\frac12-\epsilon}(N_{\min}^{12})^\sigma \|\varphi_{\lambda_1,N_1}\|_{V^2_{\pm_1}}\|\phi_{\lambda_2,N_2}\|_{V^2_{\pm_2}},
\end{align}
	and slightly bigger bound as
	\begin{align}\label{est-hh2}
		\|P_{\lambda_0}H_{N_0}(\varphi_{\lambda_1,N_1}^\dagger\beta\,\phi_{\lambda_2,N_2})	\|_{L^2_{t,x}} & \lesssim \lambda_0\left(\frac{\lambda_0}{\lambda_1}\right)^{\frac12-2\epsilon}(N_0)^\sigma \|\varphi_{\lambda_1,N_1}\|_{V^2_{\pm_1}}\|\phi_{\lambda_2,N_2}\|_{V^2_{\pm_2}},
	\end{align}
where $N_{\min}^{12}$ is the minimum of $N_1$ and $N_2$. In the region $\lambda_{\min}^{12}\ll\lambda_{\max}^{12}$, we have
\begin{align}\label{est-lh}
	\|P_{\lambda_0}H_{N_0}(\varphi_{\lambda_1,N_1}^\dagger\beta\,\phi_{\lambda_2,N_2})	\|_{L^2_{t,x}} & \lesssim (\lambda_{\min}^{12})^{\frac34-\epsilon}(\lambda_{\max}^{12})^{\frac14+\epsilon}(N_{\min})^\sigma \|\varphi_{\lambda_1,N_1}\|_{V^2_{\pm_1}}\|\phi_{\lambda_2,N_2}\|_{V^2_{\pm_2}}.
\end{align}
\end{thm}


First we decompose the modulation as follows:
\begin{align*}
P_{\lambda_0}H_{N_0}(\varphi^\dagger_{\lambda_1,N_1}\beta\phi_{\lambda_2,N_2}) & = \sum_{d\in2^{\mathbb Z}}\big( C_dP_{\lambda_0}H_{N_0}(C^{\pm_1}_{\le d}\varphi_{\lambda_1,N_1})^\dagger\beta (C^{\pm_2}_{\le d}\phi_{\lambda_2,N_2} )\\
&\qquad\qquad+ C_{<d}P_{\lambda_0}H_{N_0}(C^{\pm_1}_d\varphi_{\lambda_1,N_1})^\dagger\beta (C_{<d}^{\pm_2}\phi_{\lambda_2,N_2} )\\
&\qquad\qquad + C_{<d}P_{\lambda_0}H_{N_0}(C^{\pm_1}_{<d}\varphi_{\lambda_1,N_1})^\dagger\beta (C^{\pm_2}_d\phi_{\lambda_2,N_2})  \big)	\\
& =: \sum_{d\in2^{\mathbb Z}}\mathcal A_0+\mathcal A_1+\mathcal A_2.
\end{align*}
We consider the case $\lambda_0\ll\lambda_1\approx\lambda_2$ with $d\lesssim\lambda_0$. By Lemma \ref{resonance}, we must have $\pm_1=\pm_2$. We also note that the range of the modulation $d$ is restricted to the region $\lambda_0^{-1}\lesssim d\lesssim\lambda_0$. We begin with the $\mathcal A_0$ term.
We use the almost orthogonal decomposition by angular sectors and cubes as follows:
\begin{align*}
\|\mathcal A_0\|_{L^2_{t,x}} & \lesssim \sum_{\substack{\kappa,\kappa'\in\mathcal C_\theta \\ |\kappa-\kappa'|\lesssim\theta}}\sum_{\substack{q,q'\in\mathcal Q_{\lambda_0} \\ |q-q'|\lesssim\lambda_0}} \Big\|C_dP_{\lambda_0}H_{N_0}(C^{\pm_1}_{\le d}R_\kappa P_q\varphi_{\lambda_1,N_1})^\dagger\beta(C^{\pm_2}_{\le d}R_{\kappa'}P_{q'}\phi_{\lambda_2,N_2})\Big\|_{L^2_{t,x}},	
\end{align*}
where $\theta = \left(\frac{d\lambda_0}{\lambda_1\lambda_2} \right)^\frac12$. As we have seen \eqref{dec-di}, to exploit the null structure, we write
\begin{align*}
&\left(C^{\pm_1}_{\le d}R_\kappa P_q\varphi_{\lambda_1,N_1} \right)^\dagger\beta\left(C^{\pm_2}_{\le d}R_{\kappa'}P_{q'}\phi_{\lambda_2,N_2} \right)\\
&\qquad = \left[C^{\pm_1}_{\le d} \Big(\Pi_{\pm_1}-\Pi_{\pm_1}(\lambda_1\omega_\kappa)\Big)R_\kappa P_q\varphi_{\lambda_1,N_1}\right]^\dagger\beta(C^{\pm_2}_{\le d}R_{\kappa'}P_{q'}\phi_{\lambda_2,N_2}) \\
& \qquad\qquad +(C^{\pm_1}_{\le d}R_\kappa P_q\varphi_{\lambda_1,N_1})^\dagger\beta\left[C^{\pm_2}_{\le d}\Big(\Pi_{\pm_2}-\Pi_{\pm_2}(\lambda_2\omega_{\kappa'})\Big)R_{\kappa'}P_{q'}\phi_{\lambda_2,N_2}\right] \\
&\qquad\qquad + (C^{\pm_1}_{\le d}R_\kappa P_q\varphi_{\lambda_1,N_1})^\dagger\Pi_{\pm_1}(\lambda_1\omega_\kappa)\beta\Pi_{\pm_2}(\lambda_2\omega_{\kappa'})(C^{\pm_2}_{\le d}R_{\kappa'}P_{q'}\phi_{\lambda_2,N_2}).	
\end{align*}
Thus by Lemma \ref{lem-null}, \eqref{dirac-null} and H\"{o}lder's inequality, we gain the angle $\theta$ and then apply in order angular concentration estimates Lemma \ref{ang-con} on the lowest angular frequency term and square-summation-version of $L^4$-Strichartz estimates \eqref{cube-lin}.
\begin{align*}
\|\mathcal A_0\|_{L^2_{t,x}} 
& \lesssim 	\sum_{\substack{\kappa,\kappa'\in\mathcal C_\theta \\ |\kappa-\kappa'|\lesssim\theta}}\sum_{\substack{q,q'\in\mathcal Q_{\lambda_0} \\ |q-q'|\lesssim\lambda_0}}\theta\|C^{\pm_1}_{\le d}R_\kappa P_q\varphi_{\lambda_1,N_1}\|_{L_{t,x}^4} \|C^{\pm_2}_{\le d}R_{\kappa'}P_{q'}\phi_{\lambda_2,N_2}\|_{L^4_{t,x}} \\
& \lesssim \theta^{1-\epsilon}\left(\frac{\lambda_0}{\lambda_1}\right)^{-\epsilon}(\lambda_0\lambda_1)^\frac12 (\theta N_{\min}^{12})^\sigma \|\varphi_{\lambda_1,N_1}\|_{V^2_{\pm_1}}\|\phi_{\lambda_2,N_2}\|_{V^2_{\pm_2}} \\
& \lesssim \theta^{\sigma-\epsilon}\left(\frac{\lambda_0}{\lambda_1}\right)^{\frac12-\epsilon}(d\lambda_0)^\frac12 (N_{\min}^{12})^\sigma \|\varphi_{\lambda_1,N_1}\|_{V^2_{\pm_1}}\|\phi_{\lambda_2,N_2}\|_{V^2_{\pm_2}}.
\end{align*}
We put $\sigma>\epsilon$. (Since $\epsilon>0$ by Lemma \ref{square-sum} can be chosen arbitrarily small, $\sigma>\epsilon$ obviously implies $\sigma>0$.)  Note that $\mathcal A_1$ and $\mathcal A_2$ can be treated in the identical manner. The summation with respect to $d\lesssim \lambda_0$ gives \eqref{est-hh1}.

On the other hand, in order to get \eqref{est-hh2} we use $L^2$-duality.  Since $\pm_1  = \pm_2$, $\lambda_0 \lesssim \lambda_1 \approx \lambda_2$, and $\xi_0 = -\xi_1 + \xi_2$, Lemma \ref{resonance} leads us that
$$
\lambda_0 \gtrsim d \gtrsim \mathcal M_{\pm_1,\pm_2}(\xi_1,\xi_2) \gtrsim   \max\left( \frac{|\xi_0||\xi_1|}{\langle \xi_1 \rangle + \langle \xi_2 \rangle}\angle(\xi_0,-\xi_1)^2 , \frac{|\xi_0||\xi_2|}{\langle \xi_1 \rangle + \langle \xi_2 \rangle}\angle(\xi_0,\xi_2)^2 \right)\gtrsim\lambda_0.
$$
where $\xi_j$'s are frequencies of $\psi_j$. Then we get
\begin{align*}
&\|\mathcal A_0\|_{L^2_{t,x}}  \lesssim \sup_{\|\psi\|_{L^2_{t,x}}\lesssim1}  \sum_{\substack{\kappa_0\in\mathcal C_{\theta^*} \\ |\kappa_0\pm_2\kappa_2|\lesssim{\theta^*}}}   \sum_{\substack{\kappa_1,\kappa_2\in\mathcal C_\theta \\ |\kappa_1-\kappa_2|\approx\theta}}\sum_{\substack{q,q'\in\mathcal Q_{\lambda_0} \\ |q-q'|\lesssim\lambda_0}}\left|\mathbf I_1\right|,	
\end{align*}
where $\theta=\left(\frac{d\lambda_0}{\lambda_1\lambda_2} \right)^\frac12$, ${\theta^*}=\left(\frac{d}{\lambda_0} \right)^\frac12$, and
$$
\mathbf I_1 := \int R_{\kappa_0}C_{d}\psi_{\lambda_0,N_0}(R_{\kappa_1} P_qC^{\pm_1}_{\le d}\varphi_{\lambda_1,N_1})^\dagger\beta(R_{\kappa_2}P_{q'}C^{\pm_2}_{\le d}\phi_{\lambda_2,N_2})\,dtdx.
$$
Then as previous argument, we exploit null structure, H\"{o}lder's inequality and use angular concentration estimates on the $\psi_{\lambda_0,N_0}$ term and $L^4$-Strichartz estimates to obtain
\begin{align*}
\|\mathcal A_0\|_{L^2_{t,x}} & \lesssim \theta^{1-\epsilon}{\left(\theta^* \right)}^\sigma \left(\frac{\lambda_0}{\lambda_1}\right)^{-\epsilon}(\lambda_0\lambda_1)^\frac12 N_0^\sigma\|\varphi_{\lambda_1,N_1}\|_{V^2_{\pm_1}}\|\phi_{\lambda_2,N_2}\|_{V^2_{\pm_2}}.
\end{align*}
As \eqref{est-hh1}, the identical argument is applied to treat $\mathcal A_1,\mathcal A_2$ and hence the summation with respect to the $d\lesssim\lambda_0$ gives \eqref{est-hh2}.

In the Low$\times$High interaction, by symmetry, it suffices to deal with $\lambda_1\ll\lambda_2\approx\lambda_0$ with $d\lesssim\lambda_1$. As the previous estimate, we have $\lambda_1^{-1}\lesssim d\lesssim\lambda_1$. In the Low$\times$High interaction regime, the output frequency $\lambda_0$ is high, i.e., $\lambda_0\approx\lambda_{\max}$. To exploit the almost orthogonality by smaller cubes, we make use of $L^2$-duality as follows:
\begin{align*}
\|\mathcal A_0\|_{L^2_{t,x}} & = \sup_{\|\psi\|_{L^2_{t,x}}\lesssim1}\left|\int C_d\psi_{\lambda_0,N_0}(C^{\pm_1}_{\le d}\varphi_{\lambda_1,N_1})^\dagger\beta (C^{\pm_2}_{\le d}\phi_{\lambda_2,N_2})\,dtdx\right| \\
& \lesssim \sup_{\|\psi\|_{L^2_{t,x}}\lesssim1}\sum_{\substack{q',q''\in\mathcal Q_{\lambda_1} \\ |q'-q''|\lesssim\lambda_1}}\sum_{\substack{\kappa_0,\kappa_1,\kappa_2\in\mathcal C_{\theta^*} \\ |\kappa_1\mp_2\kappa_2|,|\kappa_0\pm_2\kappa_2|\lesssim{\theta^*}}}\left|\mathbf I_2\right|,
\end{align*}
where $\theta^* = \left( \frac{d}{\lambda_0}\right)^\frac12$ and
$$
\mathbf I_2 := \int C_dP_{q''}R_{\kappa_0}\psi_{\lambda_0,N_0}(C^{\pm_1}_{\le d}R_{\kappa_1}\varphi_{\lambda_1,N_1})^\dagger\beta(C^{\pm_2}_{\le d}P_{q'}R_{\kappa_2}\phi_{\lambda_2,N_2})\,dtdx.
$$

  Now the remainder step is very similar as the proof of \eqref{est-hh2}. Indeed, we apply in order H\"{o}lder's inequality and Lemma \ref{ang-con} and then square-sum estimates \eqref{cube-lin}.
\begin{align*}
\|\mathcal A_0\|_{L^2_{t,x}} & \lesssim 	 \sup_{\|\psi\|_{L^2_{t,x}}\lesssim1}\sum_{\substack{q',q''\in\mathcal Q_{\lambda_1} \\ |q'-q''|\lesssim\lambda_1}}\sum_{\substack{\kappa_0,\kappa_1,\kappa_2\in\mathcal C_{\theta^*} \\ |\kappa_1\mp_2\kappa_2|,|\kappa_0\pm_2\kappa_2|\lesssim{\theta^*}}}{\theta^*} \|C_dP_{q''}R_{\kappa_0}\psi_{\lambda_0,N_0}\|_{L^2_{t,x}} \\
&\qquad\qquad\qquad\qquad\qquad\times \|C^{\pm_1}_{\le d}R_{\kappa_1} \varphi_{\lambda_1,N_1}\|_{L^4_{t,x}}\|C^{\pm_2}_{\le d}P_{q'}R_{\kappa_2}\phi_{\lambda_2,N_2}\|_{L^4_{t,x}} \\
& \lesssim \left({\theta^*}\right)^{1-\epsilon}({\theta^*} N_{\min})^\sigma \left(\frac{\lambda_1}{\lambda_2}\right)^{-\epsilon}(\lambda_1\lambda_2)^\frac14\lambda_1^\frac12 \|\varphi_{\lambda_1,N_1}\|_{V^2_{\pm_1}}\|\phi_{\lambda_2,N_2}\|_{V^2_{\pm_2}} \\
& \lesssim \left({\theta^*}\right)^{\sigma-\epsilon}d^\frac12 \lambda_1^{\frac14-\epsilon}\lambda_2^{\frac14+\epsilon} (N_{\min})^\sigma \|\varphi_{\lambda_1,N_1}\|_{V^2_{\pm_1}}\|\phi_{\lambda_2,N_2}\|_{V^2_{\pm_2}},
\end{align*}
Again, we put $\sigma>\epsilon$. As the High$\times$High regime, the estimate of the $\mathcal A_1$ and $\mathcal A_2$ terms is followed in the similar way. Moreover, the summation on $d\lesssim\lambda_1$ gives \eqref{est-lh}.
This completes the proof of Proposition \ref{bi-est} in low-modulation regime.

\subsection{High modulation - Low frequency I}
The aim of this section is to prove Proposition \ref{bi-est} in the regime: $\lambda_{\min}\ll d\ll\lambda_{\max}$. It suffices to show the following bilinear estimates:
\begin{thm}
	For any $\epsilon>0$, we have
	\begin{align}\label{est-hh-mod}
	\|P_{\lambda_0}H_{N_0}(\varphi_{\lambda_1,N_1}^\dagger\beta\,\phi_{\lambda_2,N_2})\|_{L^2_{t,x}} & \lesssim \lambda_0\left(\frac{\lambda_0}{\lambda_1}\right)^{\frac14-\epsilon}	\|\varphi_{\lambda_1,N_1}\|_{V^2_{\pm_1}}\|\phi_{\lambda_2,N_2}\|_{V^2_{\pm_2}}
	\end{align}
in the High$\times$High interaction, and
\begin{align}\label{est-lh-mod}
	\|P_{\lambda_0}H_{N_0}(\varphi_{\lambda_1,N_1}^\dagger\beta\,\phi_{\lambda_2,N_2})\|_{L^2_{t,x}} & \lesssim (\lambda_{\min})^{\frac34-\epsilon}(\lambda_{\max})^{\frac14+\epsilon}\|\varphi_{\lambda_1,N_1}\|_{V^2_{\pm_1}}\|\phi_{\lambda_2,N_2}\|_{V^2_{\pm_2}}
\end{align}
in the Low$\times$High interaction.
\end{thm}


We start with the High$\times$High interaction. For $\lambda_0\ll d\ll\lambda_1$, since $\mathcal M_{\pm_1,\pm_2}\lesssim d\ll\lambda_1$, we have $\pm_1=\pm_2$ and $\mathcal M_{\pm_1,\pm_2}\lesssim\lambda_0$. Thus the angle between the support of $\widehat{\varphi}$ and $\widehat{\phi}$ is less than $(\frac{\lambda_0}{\lambda_1})^\frac12$. We first consider the $\mathcal A_0$ term. We decompose it into the following:
\begin{align}\label{mod-dec0}
\begin{aligned}
\|\mathcal A_0\|_{L^2_{t,x}} & \lesssim \|C_dP_{\lambda_0}H_{N_0}(C^{\pm_1}_{\approx d}\varphi_{\lambda_1,N_1})^\dagger\beta(C^{\pm_2}_{\ll d}\phi_{\lambda_2,N_2})\|_{L^2_{t,x}}  \\
&\qquad\qquad + \|C_dP_{\lambda_0}H_{N_0}(C^{\pm_1}_{\ll d}\varphi_{\lambda_1,N_1})^\dagger\beta(C^{\pm_2}_{\approx d}\phi_{\lambda_2,N_2})\|_{L^2_{t,x}} \\
&=: A_{0,1}+A_{0,2}. 	
\end{aligned}
\end{align}
We use the $L^2$-duality and null-form-type bound to gain $\left(\frac{\lambda_0}{\lambda_1}\right)^\frac12$. Then the almost orthogonal decomposition by smaller cubes, H\"{o}lder's inequality, Bernstein's inequality for $\psi_{\lambda_0,N_0}$, using \eqref{bdd-high-mod} for $\varphi_{\lambda_1,N_1}$ and $L^4$-Strichartz estimates \eqref{cube-lin} for $\phi_{\lambda_2,N_2}$ give us the desired estimates as follows.
 \begin{align*}
 A_{0,1} & \lesssim \sup_{\|\psi\|_{L^2_{t,x}}\lesssim1}\left(\frac{\lambda_0}{\lambda_1}\right)^\frac12\left|\int C_d\psi_{\lambda_0,N_0}(C^{\pm_1}_{\approx d}\varphi_{\lambda_1,N_1})^\dagger\beta(C^{\pm_2}_{\ll d}\phi_{\lambda_2,N_2})\,dtdx\right| \\
 & \lesssim 	\sup_{\|\psi\|_{L^2_{t,x}}\lesssim1}\left(\frac{\lambda_0}{\lambda_1}\right)^\frac12\sum_{\substack{q,q'\in\mathcal Q_{\lambda_0} \\ |q-q'|\lesssim\lambda_0}} \left|\int C_d\psi_{\lambda_0,N_0}(P_qC^{\pm_1}_{\approx d}\varphi_{\lambda_1,N_1})^\dagger\beta(P_{q'}C^{\pm_2}_{\ll d}\phi_{\lambda_2,N_2})\,dtdx\right| \\
 & \lesssim \sup_{\|\psi\|_{L^2_{t,x}}\lesssim1}\left(\frac{\lambda_0}{\lambda_1}\right)^\frac12 \sum_{\substack{q,q'\in\mathcal Q_{\lambda_0} \\ |q-q'|\lesssim\lambda_0}} \|C_d\psi_{\lambda_0,N_0}\|_{L^4_{t,x}}\|P_qC^{\pm_1}_{\approx d}\varphi_{\lambda_1,N_1}\|_{L^2_{t,x}}\|P_{q'}C^{\pm_2}_{\ll d}\phi_{\lambda_2,N_2}\|_{L^4_{t,x}} \\
 & \lesssim \left(\frac{\lambda_0}{\lambda_1}\right)^{\frac34-\epsilon}\lambda_0^\frac34 d^{\frac14} d^{-\frac12}\lambda_2^\frac12 \|\varphi_{\lambda_1,N_1}\|_{V^2_{\pm_1}}\|\phi_{\lambda_2,N_2}\|_{V^2_{\pm_2}} \\
 & \lesssim \lambda_0\left(\frac{\lambda_0}{d}\right)^\frac14\left(\frac{\lambda_0}{\lambda_1}\right)^{\frac14-\epsilon} \|\varphi_{\lambda_1,N_1}\|_{V^2_{\pm_1}}\|\phi_{\lambda_2,N_2}\|_{V^2_{\pm_2}}.
 \end{align*}
The $A_{0,2}$ tern can be also treated similarly. Now we turn our attention to the $\mathcal A_1$ term. As \eqref{mod-dec0} we get the following decomposition:
\begin{align}\label{mod-dec1}
\begin{aligned}
\|\mathcal A_1\|_{L^2_{t,x}} & \lesssim \|C_{\approx d}P_{\lambda_0}H_{N_0}(C^{\pm_1}_d\varphi_{\lambda_1,N_1})^\dagger\beta(C^{\pm_2}_{\ll d}\phi_{\lambda_2,N_2})\|_{L^2_{t,x}} \\
&\qquad\qquad + \|C_{\ll d}P_{\lambda_0}H_{N_0}(C^{\pm_1}_d\varphi_{\lambda_1,N_1})^\dagger\beta(C^{\pm_2}_{\approx d}\phi_{\lambda_2,N_2})\|_{L^2_{t,x}} \\
&=: A_{1,1}+A_{1,2}.
\end{aligned}	
\end{align}
The term $A_{1,1}$ is treated in the identical manner as the $\mathcal A_0$ term. The estimate of $A_{1,2}$ is very straightforward. Indeed, we use the null structure to obtain the factor $\left(\frac{\lambda_0}{\lambda_1}\right)^\frac12$ and the $L^2$-duality. Then simply using H\"{o}lder's inequality, Bernstein's inequality for $\psi_{\lambda_0,N_0}$ and the boundedness in the high-modulation region \eqref{bdd-high-mod} for $\varphi_{\lambda_1,N_1},\phi_{\lambda_2,N_2}$ give
\begin{align*}
A_{1,2} & \lesssim \left(\frac{\lambda_0}{\lambda_1}\right)^\frac12\sup_{\|\psi\|_{L^2_{t,x}}\lesssim1}\|C_{\ll d}\psi_{\lambda_0,N_0}\|_{L^\infty_{t,x}}\|C^{\pm_1}_d\varphi_{\lambda_1,N_1}\|_{L^2_{t,x}}\|C^{\pm_2}_{\approx d}\phi_{\lambda_2,N_2}\|_{L^2_{t,x}} \\
& \lesssim \left(\frac{\lambda_0}{\lambda_1}\right)^\frac12 \lambda_0^\frac32d^\frac12 d^{-1} \|\varphi_{\lambda_1,N_1}\|_{V^2_{\pm_1}}\|\phi_{\lambda_2,N_2}\|_{V^2_{\pm_2}}.
\end{align*}
Since we have two high-input frequencies, the estimate of $\mathcal A_2$ is exactly same as the $\mathcal A_1$ term. We omit the details. Finally, combining the bound of $\mathcal A_j$, $j=0,1,2$ and summation by $\lambda_{\min}\ll d\ll\lambda_{\max}$ gives \eqref{est-hh-mod}.

Now we are concerned with the Low$\times$High interaction. By symmetry, it is enough to consider $\lambda_1\ll\lambda_0\approx\lambda_2$ with $\lambda_1\ll d\ll \lambda_0$. The following argument is very similar as the High$\times$High interaction. In fact, we use the decomposition as \eqref{mod-dec0} and \eqref{mod-dec1}. For the Low$\times$High regime, we do not use null structure. Instead, we apply the orthogonal decomposition by cubes $q',q''\in\mathcal Q_{\lambda_1}$ to obtain
\begin{align*}
A_{0,1} 	& \lesssim \sup_{\|\psi\|_{L^2_{t,x}}\lesssim1}\left|\int C_d\psi_{\lambda_0,N_0}(C^{\pm_1}_{\approx d}\varphi_{\lambda_1,N_1})^\dagger\beta(C^{\pm_2}_{\ll d}\phi_{\lambda_2,N_2})\,dtdx\right| \\
& \lesssim \sup_{\|\psi\|_{L^2_{t,x}}\lesssim1}\sum_{\substack{q',q''\in\mathcal Q_{\lambda_1} \\ |q'-q''|\lesssim\lambda_1}}\left|\int P_{q''}C_d\psi_{\lambda_0,N_0}(C^{\pm_1}_{\approx d}\varphi_{\lambda_1,N_1})^\dagger\beta(P_{q'}C^{\pm_2}_{\ll d}\phi_{\lambda_2,N_2})\,dtdx\right|.
\end{align*}
Then H\"{o}lder's inequality, Bernstein's inequality for $\phi_{\lambda_2,N_2}$ to get $L^4$ norm and square summation by cubes followed by a simple bound for high-modulation \eqref{bdd-high-mod} give the required estimates as follows.
\begin{align*}
A_{0,1}& \lesssim \sup_{\|\psi\|_{L^2_{t,x}}\lesssim1}\sum_{\substack{q',q''\in\mathcal Q_{\lambda_1} \\ |q'-q''|\lesssim\lambda_1}} \|P_{q''}C_d\psi_{\lambda_0,N_0}\|_{L^2_{t,x}}\|C^{\pm_1}_{\approx d}\varphi_{\lambda_1,N_1}\|_{L^2_{t,x}}\|P_{q'}C^{\pm_2}_{\ll d}\phi_{\lambda_2,N_2}\|_{L^\infty_{t,x}} \\
& \lesssim \sup_{\|\psi\|_{L^2_{t,x}}\lesssim1}\sum_{\substack{q',q''\in\mathcal Q_{\lambda_1} \\ |q'-q''|\lesssim\lambda_1}} \|P_{q''}C_d\psi_{\lambda_0,N_0}\|_{L^2_{t,x}}\|C^{\pm_1}_{\approx d}\varphi_{\lambda_1,N_1}\|_{L^2_{t,x}}\lambda_1^\frac34d^\frac14\|P_{q'}C^{\pm_2}_{\ll d}\phi_{\lambda_2,N_2}\|_{L^4_{t,x}} \\
& \lesssim d^{-\frac12}\lambda_1^\frac34d^\frac14\left(\frac{\lambda_1}{\lambda_2}\right)^{-\epsilon}(\lambda_1\lambda_2)^\frac14\|\varphi_{\lambda_1,N_1}\|_{V^2_{\pm_1}}\|\phi_{\lambda_2,N_2}\|_{V^2_{\pm_2}}.
\end{align*}
On the other hand, for $A_{0,2}$, the straightforward use of H\"{o}lder's inequality, Bernstein's inequality and then \eqref{bdd-high-mod} give us
\begin{align*}
A_{0,2} & = \|C_{d}P_{\lambda_0}H_{N_0}(C^{\pm_1}_{\ll d}\varphi_{\lambda_1,N_1})^\dagger\beta(C^{\pm_2}_{\approx d}\phi_{\lambda_2,N_2})\|_{L^2_{t,x}} \\
& \lesssim \|C^{\pm_1}_{\ll d}\varphi_{\lambda_1,N_1}\|_{L^\infty_{t,x}} \|C^{\pm_2}_{\approx d}\phi_{\lambda_2,N_2}\|_{L^2_{t,x}} \\
& \lesssim \lambda_1^\frac34d^\frac14d^{-\frac12}\|\varphi_{\lambda_1,N_1}\|_{L^4_{t,x}}\|\phi_{\lambda_2,N_2}\|_{V^2_{\pm_2}} \\
& \lesssim \lambda_1\left(\frac{\lambda_1}{d}\right)^\frac14\|\varphi_{\lambda_1,N_1}\|_{V^2_{\pm_1}}\|\phi_{\lambda_2,N_2}\|_{V^2_{\pm_2}}.	
\end{align*}
Now we consider $\mathcal A_1$. We use the decomposition \eqref{mod-dec1}.
Then the $A_{1,1}$ can be treated in the similar way as $ A_{0,1}$. The estimate of $A_{1,2}$ also follows the routine of the $A_{0,1}$ term. In fact, after the $L^2$-duality we apply the orthogonal decomposition by cubes and Bernstein's inequality and then \eqref{bdd-high-mod} as follows.
\begin{align*}
A_{1,2} & \lesssim \|C_{\ll d}P_{\lambda_0,N_0}(C^{\pm_1}_{d}\varphi_{\lambda_1,N_1})^\dagger\beta(C^{\pm_2}_{\approx d}\phi_{\lambda_2,N_2})\|_{L^2_{t,x}} \\
& \lesssim \sup_{\|\psi\|_{L^2_{t,x}}\lesssim1}\sum_{\substack{q',q''\in\mathcal Q_{\lambda_1} \\ |q'-q''|\lesssim\lambda_1}}\|P_{q''}C_{\ll d}\psi_{\lambda_0,N_0}\|_{L^\infty_{t,x}}\|C^{\pm_1}_d\varphi_{\lambda_1,N_1}\|_{L^2_{t,x}}\|P_{q'}C^{\pm_2}_{\approx d}\phi_{\lambda_2,N_2}\|_{L^2_{t,x}} \\
& \lesssim \sup_{\|\psi\|_{L^2_tL^2_x}\lesssim1}\sum_{\substack{q',q''\in\mathcal Q_{\lambda_1} \\ |q'-q''|\lesssim\lambda_1}}\lambda_1^\frac32d^{\frac12}\|P_{q''}C_{\ll d}\psi_{\lambda_0,N_0}\|_{L^2_{t,x}}d^{-\frac12}\|\varphi_{\lambda_1,N_1}\|_{V^2_{\pm_1}}\|P_{q'}C^{\pm_2}_{\approx d}\phi_{\lambda_2,N_2}\|_{L^2_{t,x}} \\
& \lesssim \lambda_1^{\frac32}d^{-\frac12}\|\varphi_{\lambda_1,N_1}\|_{V^2_{\pm_1}}\|\phi_{\lambda_2,N_2}\|_{V^2_{\pm_2}}.
\end{align*}
The $\mathcal A_2$ term can be treated similarly. In fact, we can decompose it into $\mathcal A_2=A_{2,1}+A_{2,2}$ as \eqref{mod-dec0} and \eqref{mod-dec1}. Then one can deal with $A_{2,1}$ as $A_{0,2}$ and $A_{2,2}$ as $A_{1,2}$, respectively. This completes the proof of  Proposition \ref{bi-est} in the case: $\lambda_{\min}\ll d\ll\lambda_{\max}$.

\section{Non-scattering: Proof of Theorem \ref{nonscatter-thm}}

For the proof of Theorem \ref{nonscatter-thm} we have only to consider the scattering as $t \to +\infty$. We proceed by contradiction, assuming that
$$
\|\psi_{\infty}^{\ell}(0)\|_{L_x^2} > 0.
$$

Let us define functional $H(t)$ by
\begin{align*}
	H(t) =  {\rm Im} \left< \psi, \psi_{\infty}^{\ell} \right>_{L_x^2}.
\end{align*}
It is clear that $H(t)$ is uniformly bounded from the mass conservation \eqref{m-conserv}. 
Now by taking the time derivative (this can be done by a standard approximation with smooth functions) and using the self-adjointness of Dirac operator we have
\begin{align*}
	&	\frac{d}{dt}H(t) = \frac1{4\pi} {\rm Re} \left<  |\cdot |^{-1}*(\left<\psi, \beta \psi\right>_{\mathbb C^4} \beta \psi, \psi_{\infty}^{\ell} \right>_{L_x^2}.
\end{align*}
Integrating over $[t_*, t^*]$, we get
$$
H(t^*) - H(t_*) =  \int_{t_*}^{t^*}Y(t)\,dt,
$$
where $Y(t) = \frac1{4\pi} {\rm Re} \left<  |\cdot |^{-1}*\left<\psi, \beta \psi\right>_{\mathbb C^4} \beta \psi,  \psi_{\infty}^{\ell} \right>_{L_x^2}$.

We will show that
\begin{align}\label{lowerbound}
	|Y(t)| \gtrsim t^{-1}
\end{align}
for sufficiently large $t$ when $\|\psi_{\infty}^{\ell}\|_{L_x^2} > 0$.
Once \eqref{lowerbound} has been shown, due to the fixed sign of $Y$ for large time, \eqref{lowerbound} would lead us to a contradiction to the fact that $H(t)$ is uniformly bounded on time and hence complete the proof of Theorem \ref{nonscatter-thm}.

From now on we focus on the proof of \eqref{lowerbound}. To do so we reconstitute $Y$ as follows:
\begin{align*}
	Y&:= Y_1 + Y_2 + Y_3,\\	
	Y_1 &= \frac1{4\pi}  \mathcal I(\psi_{\infty}^{\ell}, \beta\psi_{\infty}^{\ell}),\\
	Y_2 &= \frac1{4\pi}  {\rm Re} \left<   |\cdot|^{-1}*\left(\left<\psi, \beta \psi\right>_{\mathbb C^4} - \left<{\psi_{\infty}^{\ell}} \beta \psi_{\infty}^{\ell}\right>_{\mathbb C^4}\right) \beta \psi_{\infty}^{\ell}, \psi_{\infty}^{\ell} \right>_{L_x^2},\\
	Y_3 &= \frac1{4\pi}  {\rm Re} \left<   |\cdot|^{-1}*\left<\psi, \beta \psi\right>_{\mathbb C^4} \beta (\psi -\psi_{\infty}^{\ell}), \psi_{\infty}^{\ell} \right>_{L_x^2}.
\end{align*}

We first deal with the $Y_1$.  By the assumption \eqref{nonscatter-condi} we have that for any $t > t_*$
\begin{align}
	\begin{aligned}\label{y1-esti}
		|Y_1| &= \frac1{4\pi} |\mathcal I(\psi_{\infty}^{\ell}, \beta\psi_{\infty}^{\ell})(t)| \ge c \frac1{4\pi}\mathcal I(\psi_{\infty}^{\ell}, \psi_{\infty}^{\ell})\\
		&= c\frac1{4\pi} \int\Big( |\cdot|^{-1}* |\psi_{\infty}^{\ell}|^2\Big)(t, x)|\psi_{\infty}^{\ell}(t, x)|^2dx\\
		&\ge c\frac1{4\pi}  (4At)^{-1}   \left(\int   \rho\left(\frac{x}{At}\right) \left|\psi_{\infty}^{\ell}(t, x)\right|^2 dx\right)^2.
	\end{aligned}
\end{align}
Here $\rho$ is the same cut-off function as previously. Let us set $\psi_{\infty}^{\ell}(0) = \varphi_{\infty}$ and for the sake of simplicity, let us denote $\Pi_{\pm}\varphi_{\infty}$ by $\varphi_\pm$. Then  $\psi_{\infty}^{\ell} = e^{-it\left< \nabla\right>} \varphi_+ - e^{it\left< \nabla\right>} \varphi_-$ and we obtain
\begin{align*}
	\int \rho\left(\frac{x}{At}\right)   \left|\psi_{\infty}^{\ell}(x)\right|^2 dx&= \int \rho\left(\frac{x}{At}\right) |e^{-it\left< \nabla\right>} \varphi_+ - e^{it\left< \nabla\right>} \varphi_- |^2dx\\
	&= \int   \rho\left(\frac{x}{At}\right)  \left(|e^{-it\left< \nabla\right>} \varphi_+|^2  +    |e^{it\left< \nabla\right>} \varphi_- |^2  -2  {\rm Re} \left<e^{-it\left< \nabla\right>} \varphi_+, e^{it\left< \nabla\right>} \varphi_- \right>_{\mathbb C^4} \right) dx.
\end{align*}
We handle the last term in the integrand as follows:
\begin{align*}
	\int \rho\left(\frac{x}{At} \right) \left<e^{-it\left< \nabla\right>} \varphi_+, e^{it\left< \nabla\right>} \varphi_- \right>_{\mathbb C^4}  dx &=  \left< \rho\left(\frac{x}{At} \right)e^{-it\left< \nabla\right>} \varphi_+, e^{it\left< \nabla\right>}\Pi_{-} \varphi_{\infty}\right>_{L_x^2}\\
	&=  \left<\Pi_{-} \left(\rho\left(\frac{x}{At} \right)  e^{-it\left< \nabla\right>} \varphi_+ \right), e^{it\left< \nabla\right>} \varphi_\infty\right>_{L_x^2}\\
	&=  \left<\rho\left(\frac{x}{At} \right)e^{-it\left< \nabla\right>} \Pi_{-} \Pi_+\varphi_\infty,  e^{it\left< \nabla\right>} \varphi_\infty\right>_{L_x^2}\\
	&\qquad\;\; + \left<  \left[\Pi_{-}, \rho\left(\frac{x}{At} \right)\right]  e^{-it\left< \nabla\right>} \varphi_+, e^{it\left< \nabla\right>} \varphi_\infty\right>_{L_x^2}\\
	&= - \frac12  \left<\left[\frac{\alpha^j\partial_j - \beta}{\left< \nabla\right>},\rho\left(\frac{x}{At} \right)\right]e^{-it\left< \nabla\right>} \varphi_+,   e^{it\left< \nabla\right>} \varphi_{\infty}\right>_{L_x^2}.
\end{align*}
Here $\left[A ,B \right] := AB - BA$. We used properties of \eqref{proj} for the last integral. On the other hand, Plancherel's theorem yields
\begin{align*}
	&\left|2{\rm Re}\int \left<e^{-it\left< \nabla\right>} \varphi_+ , e^{it\left< \nabla\right>}\varphi_-\right>_{\mathbb C^4} dx\right| =  \left|\int \left<\left[\frac{\alpha^j\partial_j - \beta}{\left< \nabla\right>},\rho\left(\frac{x}{At} \right)\right] \varphi_+, e^{it\left< \nabla\right>} \varphi_{\infty}\right> dx\right|\\
	&\quad\lesssim (At)^3 \iint \left| \alpha\cdot\left(\frac{\xi}{\left<\xi\right>} - \frac{\eta}{\left<\eta\right>} \right) - \beta\left(\frac1{\left<\xi\right>} - \frac1{\left<\eta\right>}\right) \right|\Big|\widehat{\rho}\left(At(\xi-\eta)\right)\Big|  \big|\Pi_+(\eta)\widehat{\varphi_{\infty}}(\eta)\big|    \left|\widehat{\varphi_{\infty}}(-\xi)\right| d\eta d\xi\\
	&\quad\lesssim (At)^3 \iint \frac{|\xi-\eta|}{\left<\xi\right>} \Big|\widehat{\rho}\left(At(\xi-\eta)\right)\Big|  \big|\Pi_+(\eta)\widehat{\varphi_{\infty}}(\eta)\big|    \left|\widehat{\varphi_{\infty}}(-\xi)\right| d\eta d\xi\\
	&\quad\lesssim (At)^{-1} \|\varphi_{\infty}\|_{L_x^2}^2.
\end{align*}
From this we can choose sufficiently large $t$ so that
\begin{align}\label{era}
	\left|2{\rm Re}\int \left<e^{-it\left< \nabla\right>} \varphi_+ , e^{it\left< \nabla\right>}\varphi_-\right>_{\mathbb C^4} dx\right| \le \frac1{10}\|\varphi_{\infty}\|_{L_x^2}.
\end{align}

Let us set $\varphi_{\kappa,r} := P_{\kappa^{-1}<\cdot\le \kappa}\left(\rho\left(\frac{\cdot}{r}\right) \varphi_{\infty}\right)$, $\phi(\xi) := x\cdot\xi + t\left<\xi\right>$, and $L(\xi) := \nabla_{\xi} \left(\frac{\nabla_{\xi}\phi}{|\nabla_{\xi}\phi|^2}\right)$. If $|x| \ge At$, by using integration by parts twice we have
\begin{align}
	\begin{aligned}\label{decay-space}
		&\left|e^{-it\left< \nabla\right>} \Pi_{+}\varphi_{\kappa,r}(x) \right| =  C \left|\int \nabla_{\xi}\phi e^{i\phi} \frac{\nabla_{\xi}\phi}{|\nabla_{\xi}\phi|^2} \widehat{\varphi_{\kappa,r}^+}(\xi)\, d \xi \right|\\
		&\quad = C  \left| \int  e^{i\phi} L(\xi) \widehat{\varphi_{\kappa,r}^+}(\xi)\, d \xi  + \int  e^{i\phi} \frac{\nabla_{\xi}\phi}{|\nabla_{\xi}\phi|^2} \nabla_{\xi} \widehat{\varphi_{\kappa,r}^+}(\xi)\, d \xi\right|\\
		&\quad \le C \left| \int  e^{i\phi} L(\xi)^2 \widehat{\varphi_{\kappa,r}^+}(\xi) + e^{i\phi} \nabla_{\xi} L(\xi) \frac{\nabla_{\xi}\phi}{|\nabla_{\xi}\phi|^2}  \widehat{\varphi_{\kappa,r}^+}(\xi) + e^{i\phi} L(\xi)^2 \nabla_{\xi}\widehat{\varphi_{\kappa,r}^+}(\xi)\, d \xi \right| \\
		&\quad\qquad +  C\left|\int  e^{i\phi} \nabla_{\xi} \left(\frac{1}{|\nabla_{\xi}\phi|^2}\right) \nabla_{\xi} \widehat{\varphi_{\kappa,r}^+}(\xi)  +   e^{i\phi}  \frac{1}{|\nabla_{\xi}\phi|^2} \nabla_{\xi}^2 \widehat{\varphi_{\kappa,r}^+}(\xi)\, d \xi\right|\\
		&\quad \le C(A,\kappa,r)|x|^{-2}\|\varphi_+\|_{L_x^2}.
	\end{aligned}
\end{align}
Since for large $\kappa$ and $r$, $\|\varphi_{\infty} - \varphi_{\kappa,r}\|_{L^2_x} \le \frac1{10}\|\varphi_{\infty}\|_{L^2_x}$, by \eqref{era} and \eqref{decay-space}, we have
\begin{align*}
	\left\|\rho\left(\frac{\cdot}{At}\right) \psi_{\infty}^{\ell}\right\|_{L_x^2}^2 &\ge \left\|\rho\left(\frac{\cdot}{At}\right)  e^{-it\left< \nabla\right>} \Pi_{+}\varphi_{\infty} \right\|_{L_x^2}^2 - \frac{1}{10}\|\varphi_{\infty}\|_{L_x^2}^2\\
	&\ge \left\|\rho\left(\frac{\cdot}{At}\right) e^{-it\left< \nabla\right>} \Pi_{+}\varphi_{\kappa,r} \right\|_{L_x^2}^2 - \frac{11}{50}\|\varphi_{\infty}\|_{L_x^2}^2 - \frac{1}{10}\|\varphi_{\infty}\|_{L_x^2}^2\\
	&\ge \|e^{-it\left< \nabla\right>} \Pi_{+} \varphi_{\kappa,r}\|_{L_x^2}^2 - \int \left( 1- \rho\left(\frac{\cdot}{At}\right)  \right) |e^{-it\left< \nabla\right>}\Pi_{+} \varphi_{\kappa,r}|^2 dx - \frac{8}{25}\|\varphi_{\infty}\|_{L_x^2}^2\\
	&\ge \frac{17}{25}\|\varphi_{\infty}\|_{L_x^2}^2 -  C(A,r)\|\varphi_{\infty}\|_{L_x^2}^2\int_{|x|\ge At} |x|^{-4}dx\\
	&\ge \frac35\|\varphi_{\infty}\|_{L_x^2}^2.
\end{align*}
This together with \eqref{y1-esti} leads us to
\begin{align}\label{y1}
	|Y_1|  \gtrsim |t|^{-1} .
\end{align}

Now let us turn to $Y_2,\;Y_3$. To treat them we need a time decay estimate for the linear solutions and $L_x^\infty$ estimates for the potential term.

\begin{lem}[see Lemma 4.2 of \cite{chooz}]\label{time-decay}
	Let $f \in B_{1,1}^\frac52$. Then
	\begin{align*}
		\|e^{\pm it\left<\nabla\right>} f\|_{L_x^\infty} \lesssim t^{-\frac32} \|f \|_{B_{1,1}^\frac52}.
	\end{align*}
\end{lem}
Here $B_{1,1}^\frac52$ is the inhomogeneous Besov space defined by $\left\{ f : \|f\|_{B_{1,1}^\frac52}:= \sum_{\lambda\ge 1}\lambda^\frac52 \|P_{\lambda} f\|_{L_x^1}  < \infty\right\}$.

\begin{lem}\label{infty-esti-1}
	For any $\mathbb C$-valued functions $u \in L_x^2 \cap L_x^\infty$ we have
	\begin{align*}
		\||x|^{-1} *|u|^2 \|_{L_x^\infty} \lesssim \|u\|_{L_x^2} \|u\|_{L_x^6}.
	\end{align*}
\end{lem}
This lemma can be readily shown by a standard optimization.

If $u_1 \neq u_2$, then we need a regularity and a space-decay assumption.
\begin{lem}[Lemma 3.2 of \cite{chooz}]\label{infty-weight}
	Let  $\mathbb C$-valued functions $u_1 \in L_x^2$ and $u_2 \in L_x^6$. Then for any $0<\epsilon<1$, we get
	\begin{align*}
		\||x|^{-1} * (|u|^2) \|_{L_x^\infty} \lesssim \|u\|_{L_x^\frac{6}{2-\varepsilon}} \|u\|_{L_x^\frac{6}{2+\varepsilon}}.
	\end{align*}
\end{lem}

From Lemmas  \ref{time-decay}, \ref{infty-weight}, and mass conservation it follows that
\begin{align*}
	|Y_2| &\lesssim   \left| \left<   V*\Big(\left<\psi, \beta \psi\right>_{\mathbb C^4} - \left<\psi_{\infty}^{\ell}, \beta \psi_{\infty}^{\ell}\right>_{\mathbb C^4}\Big) \beta \psi_{\infty}^{\ell}, \psi_{\infty}^{\ell} \right>_{L_x^2}\right|\\
	&\lesssim  \left| \int  \Big(\left<\psi, \beta \psi\right>_{\mathbb C^4} - \left<\psi_{\infty}^{\ell}, \beta \psi_{\infty}^{\ell}\right>_{\mathbb C^4}\Big) V*\left<\beta \psi_{\infty}^{\ell}, \psi_{\infty}^{\ell} \right>_{\mathbb C^4} dx\right|\\
	&\lesssim  \|\psi - \psi_{\infty}^{\ell}\|_{L_x^2} \left( \|\psi\|_{L_x^2} + \|\psi_{\infty}^{\ell}\|_{L_x^2} \right)\|\psi_{\infty}^{\ell}\|_{L_x^2}\|\psi_{\infty}^{\ell}\|_{L_x^6}\\
	&\lesssim  \|\psi - \psi_{\infty}^{\ell}\|_{L_x^2} \left( \|\psi\|_{L_x^2} + \|\psi_{\infty}^{\ell}\|_{L_x^2} \right)\|\psi_{\infty}^{\ell}\|_{L_x^2}^{\frac43}\|\psi_{\infty}^{\ell}\|_{L_x^\infty}^\frac23\\
	&\lesssim \|\psi - \psi_{\infty}^{\ell}\|_{L_x^2} \left( \|\psi_0\|_{L_x^2} + \|\varphi_{\infty}\|_{L_x^2} \right)\|\varphi_{\infty}\|_{L_x^2}^{\frac43}t^{-1}\left(\|\varphi_+\|_{B_{1,1}^\frac52} + \|\varphi_-\|_{B_{1,1}^\frac52}\right)^{\frac23}
\end{align*}
and
\begin{align*}
	|Y_3| &\lesssim \left| \left<   V*\left<\psi, \beta \psi\right>_{\mathbb C^4} \beta (\psi -\psi_{\infty}^{\ell}), \psi_{\infty}^{\ell} \right>_{L_x^2} \right|\\
	&\lesssim \left| \int   \left<\psi, \beta \psi\right>_{\mathbb C^4} \left(V* \left<\beta (\psi -\psi_{\infty}^{\ell}),\psi_{\infty}^{\ell} \right>_{\mathbb C^4}\right)  dx\right|\\
	&\lesssim \|\psi\|_{L_x^2}^2 \left\| \left|\left<\beta (\psi -\psi_{\infty}^{\ell}),\psi_{\infty}^{\ell} \right>_{\mathbb C^4} \right|^{\frac12} \right\|_{L_x^{\frac6{2-\varepsilon}}}  \left\| \left|\left<\beta (\psi -\psi_{\infty}^{\ell}),\psi_{\infty}^{\ell} \right>_{\mathbb C^4} \right|^{\frac12} \right\|_{L_x^{\frac6{2+\varepsilon}}}\\
	&\lesssim \|\psi\|_{L_x^2}^2 \left\| \left<\beta (\psi -\psi_{\infty}^{\ell}),\psi_{\infty}^{\ell} \right>_{\mathbb C^4}  \right\|_{L_x^{\frac3{2-\varepsilon}}}^{\frac12}  \left\| \left<\beta (\psi -\psi_{\infty}^{\ell}),\psi_{\infty}^{\ell} \right>_{\mathbb C^4}  \right\|_{L_x^{\frac3{2+\varepsilon}}}^{\frac12}\\
	&\lesssim \|\psi\|_{L_x^2}^2 \left\|  \psi -\psi_{\infty}^{\ell} \right\|_{L_x^2}  \|\psi_{\infty}^{\ell} \|_{L_x^{\frac6{1-\varepsilon}}}^{\frac12} \|\psi_{\infty}^{\ell} \|_{L_x^{\frac6{1+\varepsilon}}}^{\frac12}\\
	&\lesssim \|\psi\|_{L_x^2}^2 \|\psi-\psi_{\infty}^{\ell}\|_{L_x^2}\|\psi_{\infty}^{\ell}\|_{L_x^2}^\frac13 \|\psi_{\infty}^{\ell}\|_{L_x^\infty}^\frac23\\
	&\lesssim \|\psi_0\|_{L_x^2}^{2}\|\psi-\psi_{\infty}^{\ell}\|_{L_x^2}\|\varphi_{\infty}\|_{L_x^2}^{\frac13}t^{-1}\left(\|\varphi_+\|_{B_{1,1}^\frac52} + \|\varphi_-\|_{B_{1,1}^\frac52}\right)^{\frac23}
\end{align*}
Therefore we get
\begin{align}\label{yj}
	|Y_j| = o(t^{-1})
\end{align}
for $j=2,3$. Then \eqref{y1} and \eqref{yj} conclude \eqref{lowerbound}.

This completes the proof of Theorem \ref{nonscatter-thm}.

\section*{Acknowledgements}
This work was supported by NRF-2021R1I1A3A04035040(Republic of Korea).


\end{document}